\documentclass[12pt]{amsart}

\usepackage{amsmath,amsfonts,amssymb,amsthm,amscd,epsfig,comment}
\usepackage[all]{xy}

\allowdisplaybreaks

\def\C{\mathbb{C}}

\def\Z{\mathbb{Z}}

\def\id{\mathrm{id}}

\def\g{\ensuremath{\mathfrak{g}}}

\def\v{\mathbf{v}}
\def\w{\mathbf{w}}
\def\a{\mathbf{a}}

\def\i{\mathbf{i}}
\def\j{\mathbf{j}}
\def\I{\mathbf{I}}

\def\tvT{\ensuremath{\widetilde{\mathfrak{T}}(\v ; \w^1, \dots, \w^n)}}
\def\vT{\ensuremath{{\mathfrak{T}}(\v ; \w^1, \dots, \w^n)}}
\def\tT{\ensuremath{\widetilde{\mathfrak{T}}(\w^1, \dots, \w^n)}}
\def\T{\ensuremath{\mathfrak{T}(\w^1, \dots, \w^n)}}
\def\ke{{\tilde e}}
\def\kf{{\tilde f}}

\DeclareMathOperator{\im}{im} 
\DeclareMathOperator{\Hom}{Hom}

\DeclareMathOperator{\End}{End}

\DeclareMathOperator{\inc}{in}
\DeclareMathOperator{\out}{out}

\DeclareMathOperator{\tr}{tr}

\DeclareMathOperator{\wt}{wt}

\DeclareMathOperator{\Coker}{Coker}

\DeclareMathOperator{\Id}{Id}
\DeclareMathOperator{\diag}{diag}

\DeclareMathOperator{\hw}{hw}
\DeclareMathOperator{\bdim}{\mathbf{dim}}
\DeclareMathOperator{\Ob}{Ob}

\newtheorem{theo}{Theorem}[section]
\newtheorem{prop}[theo]{Proposition}
\newtheorem{lem}[theo]{Lemma}
\newtheorem{cor}[theo]{Corollary}

\newtheorem{defin}[theo]{Definition}
\newtheorem*{rem*}{Remark}

\newtheorem{example}[theo]{Example}

\numberwithin{equation}{section}

\begin{document}
\title[Crystals, quiver varieties, and coboundary categories]
{Crystals, quiver varieties, and coboundary categories for Kac-Moody
algebras}
\author{Alistair Savage}
\address{University of Ottawa\\
Ottawa, Ontario \\ Canada} \email{alistair.savage@uottawa.ca}
\date{December 1, 2008}

\begin{abstract}
Henriques and Kamnitzer have defined a commutor for the category of
crystals of a finite-dimensional complex reductive Lie algebra that
gives it the structure of a coboundary category (somewhat analogous
to a braided monoidal category). Kamnitzer and Tingley then gave an
alternative definition of the crystal commutor, using Kashiwara's
involution on Verma crystals, that generalizes to the setting of
symmetrizable Kac-Moody algebras. In the current paper, we give a
geometric interpretation of the crystal commutor using quiver
varieties. Equipped with this interpretation we show that the
commutor endows the category of crystals of a symmetrizable
Kac-Moody algebra with the structure of a coboundary category,
answering in the affirmative a question of Kamnitzer and Tingley.
\end{abstract}

\maketitle \thispagestyle{empty}

\tableofcontents

\section*{Introduction}
Let $\g$ be a symmetrizable Kac-Moody algebra and $U_q(\g)$ the
corresponding quantum group (or quantized universal enveloping
algebra).  Introduced by Kashiwara, crystals can be thought of as a
combinatorial model of representations of $U_q(\g)$ arising from the
limit as $q$ tends to zero.  To each representation of $U_q(\g)$ is
associated a crystal graph.  Roughly speaking, the crystal graph is
an edge-colored directed graph in which a certain basis (the global,
or canonical, basis) of the representation is replaced by a set of
vertices and the action of the Chevalley generators is replaced by
colored arrows. Arrows are labeled by simple roots of $\g$ and to
each vertex is associated a weight of $\g$.  One can take the tensor
product of two crystals and this operation corresponds to the tensor
product of the corresponding representations.  The vertex set of the
tensor product crystal is the Cartesian product of the two original
vertex sets. With this operation, the category of $\g$-crystals
becomes a monoidal category.  As for representations of $U_q(\g)$,
the tensor product of crystals is not symmetric.  That is, the map
$(b_1, b_2) \mapsto (b_2, b_1)$ is not a morphism of crystals in
general.

Recall that a braided monoidal category is a monoidal category
$\mathcal{C}$ equipped with a natural isomorphism $\sigma^{br}_{V,W}
: V \otimes W \to W \otimes V$ for all $V,W \in \Ob \mathcal{C}$ and
such that the diagram
\[
\xymatrix{ {U \otimes V \otimes W} \ar[rr]^{\Id \otimes
\sigma^{br}_{V,W}} \ar[drr]^{\sigma^{br}_{U \otimes V,W}} & &
{U \otimes W \otimes V} \ar[d]^{\sigma^{br}_{U,W} \otimes \Id} \\
& & {W \otimes U \otimes V}}
\]
commutes for all $U, V, W \in \Ob \mathcal{C}$.  Such a
$\sigma^{br}$ is called a \emph{braiding} and it induces an action
of the braid group on multiple tensor products. We refer the reader
to \cite[{\S}5.2]{CP95} for further details.

For $\g$ a finite-dimensional complex Lie algebra, the category of
representations of $U_q(\g)$ has a natural braiding constructed
using the universal $R$-matrix, an element of $U_h(\g) \otimes
U_h(\g)$ where $U_h(\g)$ is the formal completion of $U_q(\g)$.  The
braiding is given by the map $\text{flip} \circ R$ where
$\text{flip} : V \otimes W \to W \otimes V$ is given by $v \otimes w
\mapsto w \otimes v$.  However, this braiding does not pass to the
$q \to 0$ limit. In other words, it does not induce a braiding on
the category of crystals.  In fact, one can show that no such
braiding exists. That is, the category of $\g$-crystals cannot be
given the structure of a braided monoidal category for nontrivial
$\g$.  However, it can be given an analogous structure as we now
describe.

A \emph{coboundary} (or \emph{cactus}) \emph{category} is a monoidal
category $\mathcal{C}$ equipped with a natural isomorphism
$\sigma_{V,W} : V \otimes W \to W \otimes V$ for all $V, W \in \Ob
\mathcal{C}$, called a \emph{commutor}, such that $\sigma_{W,V}
\circ \sigma_{V,W} = \Id$ and the diagram
\begin{equation} \label{eq:cactus-relation}
\xymatrix{ U \otimes V \otimes W \ar[rr]^{\Id \otimes \sigma_{V,W}}
\ar[d]_{\sigma_{U,V} \otimes \Id} & & U \otimes W \otimes V
\ar[d]^{\sigma_{U, W \otimes V}} \\
V \otimes U \otimes W \ar[rr]^{\sigma_{V \otimes U, W}} & & W
\otimes V \otimes U}
\end{equation}
commutes for all $U, V, W \in \Ob \mathcal{C}$. The commutativity of
\eqref{eq:cactus-relation} is called the \emph{cactus relation}. A
commutor satisfying these conditions induces an action of the
\emph{cactus group} (see \cite{HK06}) on multiple tensor products.

Drinfel'd \cite{Dri90} has shown that one can use the $R$-matrix to
construct a commutor satisfying the cactus relation in the category
$U_q(\g)$ via a process he calls \emph{unitarization}.  If one
defines $R' = R(\text{flip}(R) R)^{-1/2}$ where the square root is
with respect to the $h$ filtration on $U_h(\g) \otimes U_h(\g)$,
then the map $\text{flip} \circ R'$ is a commutor.  In \cite{HK06},
Henriques and Kamnitzer, following an idea of A. Berenstein, defined
a commutor on the category of representations of $U_q(\g)$ and the
category of $\g$-crystals for $\g$ a finite-dimensional complex
reductive Lie algebra. Their definition involved the
Sch\"utzenberger involution which only exists in finite type. It was
shown by Kamnitzer and Tingley in \cite{KT07} that a particular case
of Henriques and Kamnitzer's construction agrees with Drinfel'd's
commutor and that this unitarization does pass to the $q \to 0$
limit.  That is, it induces the structure of a coboundary category
on the category of $\g$-crystals for $\g$ of finite type. In
\cite{KT06}, Kamnitzer and Tingley gave an alternative definition of
the commutor using Kashiwara's involution. The new approach has the
benefit of defining the involution $\sigma$ for arbitrary
symmetrizable Kac-Moody algebras.  Thus, there exist two
combinatorial definitions of the commutor: the definition of
\cite{HK06} where the cactus relation can easily be seen to hold but
which does not generalize to the Kac-Moody setting, and the
definition of \cite{KT06} which does generalize to the Kac-Moody
setting but for which the cactus relation has not been shown to hold
(for non-finite type). Kamnitzer and Tingley posed the natural
question of whether or not the commutor, extended to symmetrizable
Kac-Moody algebras via the second definition, satisfies the cactus
relation in the more general setting.

The goal of the current paper is twofold.  First, we give a
geometric interpretation of the crystal commutor using the quiver
varieties of Lusztig and Nakajima.  These are varieties associated
to quivers (directed graphs) constructed from the Dynkin graph of a
symmetrizable Kac-Moody algebra $\g$.  The set of irreducible
components of these varieties can be given the structure of a
$\g$-crystal in a natural geometric way.  We give a geometric
characterization of these irreducible components and use this
description to analyze the action of the crystal commutor described
above.  In doing so, we attain our second goal.  Namely, we answer
the above question of Kamnitzer and Tingley in the affirmative: the
crystal commutor satisfies the cactus relation for an arbitrary
symmetrizable Kac-Moody algebra $\g$ and therefore endows the
category of $\g$-crystals with the structure of a coboundary
category.  The key ingredient in the proof is that in the language
of quiver varieties, the two compositions of commutors appearing in
the cactus relation \eqref{eq:cactus-relation} both correspond to
taking adjoints of quiver representations (at least when we restrict
them to highest weight elements) and are therefore equal.

The paper is organized as follows.  In
Section~\ref{sec:coboundary-categories} we review the definition of
the crystal commutor using Kashiwara's involution. In
Section~\ref{sec:qv} we introduce the quiver varieties of Lusztig,
Malkin and Nakajima. The geometric realization of the crystals
corresponding to the lower half of the quantized enveloping algebra
of a symmetric Kac-Moody algebra and its integrable highest weight
representations is given in Section~\ref{sec:geometric-crystals}. In
Section~\ref{sec:geom-commutor} we discuss various characterizations
of the irreducible components of quiver varieties and examine how
the crystal commutor acts on their irreducible components. Equipped
with a precise description of this action, we prove that the
commutor satisfies the cactus relation in
Section~\ref{sec:cactus-relation}. Finally, in
Section~\ref{sec:non-simply-laced}, we extend our results to the
case of Kac-Moody algebras with symmetrizable (rather than
symmetric) Cartan matrices by a well-known ``folding'' argument.

The author would like to thank J. Kamnitzer for useful discussions
during the writing of this paper and P. Tingley for pointing out an
error in an earlier version. He would also like to thank the referee
for helpful comments and suggestions.  This research was supported
by the Natural Sciences and Engineering Research Council (NSERC) of
Canada.


\section{Crystals and coboundary categories}
\label{sec:coboundary-categories}

In this section we introduce the crystal commutor as defined by
Kamnitzer and Tingley in \cite{KT06}.  It was defined in a different
manner for the case of finite-dimensional complex reductive Lie
algebras by Henriques and Kamnitzer in \cite{HK06}.  We refer the
reader to~\cite{Sav08b} for a more detailed overview of the topic.
In the current paper, by the \emph{category of $\g$-crystals} for a
symmetrizable Kac-Moody algebra $\g$, we mean the category
consisting of those crystals $B$ such that each connected component
of $B$ is isomorphic to some $B_\lambda$, the crystal corresponding
to the irreducible highest weight $U_q(\g)$-module of highest weight
$\lambda$, where $\lambda$ is a dominant integral weight. For the
rest of this paper, the word \emph{crystal} means either an object
in this category or the crystal $B_\infty$ corresponding to the
lower half $U_q(\g)^-$ of the quantized universal enveloping
algebra.

\subsection{Kashiwara's involution}
\label{sec:kash-involution}

Let $\g$ be a symmetrizable Kac-Moody algebra and let $B_\infty$ be
the $\g$-crystal corresponding to the lower half $U_q^-(\g)$ of the
quantized universal enveloping algebra.  Let $* : U_q(\g) \to
U_q(\g)$ be the anti-automorphism given by
\begin{align*}
q^* &= q, \\
e_i^* &= e_i, \\
f_i^* &= f_i, \\
q(h)^* &= q(-h).
\end{align*}
The map $*$ sends $U_q^-(\g)$ to $U_q^-(\g)$ and induces a map $* :
B_\infty \to B_\infty$ (see \cite[\S 8.3]{Kas95}).  Setting
\begin{align*}
\ke_i^*(b) &= (\ke_i(b^*))^*, \\
\kf_i^*(b) &= (\kf_i(b^*))^*, \\
\tilde \varepsilon_i^*(b) &= \varepsilon_i(b^*), \\
\tilde \varphi_i^*(b) &= \varphi_i(b^*),
\end{align*}
gives $B_\infty$ another crystal structure.  We call the map $*$
\emph{Kashiwara's involution}.

Let $B_\lambda$ be the $\g$-crystal corresponding to the irreducible
highest weight $U_q(\g)$-module of highest weight $\lambda$ and let
$b_\lambda$ be its highest weight element. We recall the tensor
product rule for crystals.
\begin{align*}
\ke_i(b_1 \otimes b_2) &=
\begin{cases}
    \ke_i b_1 \otimes b_2 & \text{if } \varphi_i(b_1) \ge
    \varepsilon_i(b_2) \\
    b_1 \otimes \ke_i b_2 & \text{if } \varphi_i(b_1) <
    \varepsilon_i(b_2)
\end{cases},\\
\kf_i(b_1 \otimes b_2) &=
\begin{cases}
    \kf_i b_1 \otimes b_2 & \text{if } \varphi_i(b_1) >
    \varepsilon_i(b_2) \\
    b_1 \otimes \kf_i b_2 & \text{if } \varphi_i(b_1) \le
    \varepsilon_i(b_2)
\end{cases},\\
\wt (b_1 \otimes b_2) &= \wt(b_1) + \wt(b_2),\\
\varepsilon_i(b_1 \otimes b_2) &= \max (\varepsilon_i(b_1),
\varepsilon_i(b_2) - \left< h_i, \wt (b_1) \right> ),\\
\varphi_i(b_1 \otimes b_2) &= \max (\varphi_i(b_2), \varphi_i(b_1) +
\left< h_i, \wt (b_2) \right> ).
\end{align*}
For two dominant integral weights $\lambda$ and $\mu$, there is an
inclusion of crystals $B_{\lambda + \mu} \hookrightarrow B_\lambda
\otimes B_\mu$ sending $b_{\lambda + \mu}$ to $b_\lambda \otimes
b_\mu$.  It follows from the tensor product rule that the image of
this inclusion contains all elements of the form $b \otimes b_\mu$
for $b \in B_\lambda$.  Thus we define a map
\[
\iota^{\lambda+\mu}_\lambda : B_\lambda \to B_{\lambda + \mu}
\]
which sends $b \in B_\lambda$ to the inverse image of $b \otimes
b_\mu$ under the inclusion $B_{\lambda + \mu} \hookrightarrow
B_\lambda \otimes B_\mu$.  While this map is not a morphism of
crystals, it takes $b_\lambda$ to $b_{\lambda + \mu}$ and is
$\ke_i$-equivariant.  Here \emph{$\ke_i$-equivariant} means that
$\iota^{\lambda+\mu}_\lambda (\ke_i b) = \ke_i
\iota^{\lambda+\mu}_\lambda (b)$ for all $i$ (note that it follows
that $\ke_i b=0$ whenever $\ke_i \iota^{\lambda+\mu}_\lambda
(b)=0$). This notion of $\ke_i$-equivariant is sometimes called
\emph{$\ke_i$-strict}

The maps $\iota^{\lambda + \mu}_\lambda$ make the family of crystals
$B_\lambda$ into a directed system and the crystal $B_\infty$ can be
viewed as the limit of this system.  We have $\ke_i$-equivariant
maps $\iota^\infty_\lambda : B_\lambda \to B_\infty$ which we will
simply denote by $\iota^\infty$ when it will cause no confusion.  We
define $\iota^\infty : B \to B_\infty$ for an arbitrary $\g$-crystal
$B$ by setting $\iota^\infty(b) = \iota^\infty_\lambda(b)$ if the
connected component of $B$ containing $b$ is isomorphic to
$B_\lambda$.  This extended map $\iota^\infty$ is also
$\ke_i$-equivariant.  Define $\varepsilon^* : B_\infty \to P_+$ by
\[
\varepsilon^*(b) = \min \{\lambda\ |\ b \in
\iota^\infty(B_\lambda)\}
\]
where we put the usual order on $P_+$, the positive weight lattice
of $\g$, given by $\lambda \ge \mu$ if and only if $\lambda - \mu
\in P_+$. Recall that we also have the map $\varepsilon : B_\infty
\to P_+$ given by $\left< \alpha_i^\vee, \varepsilon(b)\right> =
\varepsilon_i(b)$. Then by \cite[Prop.~8.2]{Kas95}, Kashiwara's
involution preserves weights and satisfies
\begin{equation}
\varepsilon^*(b) = \varepsilon(b^*).
\end{equation}


\subsection{The crystal commutor}
\label{sec:crystal-commutor-def}

Consider the crystal $B_\lambda \otimes B_\mu$.  Since $\varphi(b) =
\varepsilon(b) + \wt(b)$ for all $b \in B_\lambda$, we have that
$\varphi(b_\lambda) = \wt(b_\lambda) = \lambda$. It follows from the
tensor product rule for crystals that the highest weight elements of
$B_\lambda \otimes B_\mu$ are those elements of the form $b_\lambda
\otimes b$ for $b \in B_\mu$ with $\varepsilon(b) \le \lambda$. Thus
$\varepsilon^*(b^*) = \varepsilon(b) \le \lambda$ and so, by the
definition of $\varepsilon^*$, we have $b^* \in
\iota^\infty(B_\lambda)$. So we can consider $b^*$ as an element of
$B_\lambda$.  Furthermore, $\varepsilon(b^*) = \varepsilon^*(b) \le
\mu = \varphi(b_\mu)$ since $b \in B_\mu$. Thus $b_\mu \otimes b^*$
is a highest weight element of $B_\mu \otimes B_\lambda$.  Since
$B_\lambda \otimes B_\mu \cong B_\mu \otimes B_\lambda$ as crystals,
we can make the following definition.

\begin{defin}[{\cite{KT06}}]
Let $\sigma_{B_\lambda,B_\mu} : B_\lambda \otimes B_\mu
\stackrel{\cong}{\to} B_\mu \otimes B_\lambda$ be the crystal
isomorphism given uniquely by $\sigma_{B_\lambda, B_\mu}(b_\lambda
\otimes b) = b_\mu \otimes b^*$ for $b_\lambda \otimes b$ a highest
weight element of $B_\lambda \otimes B_\mu$. The map
$\sigma_{B_\lambda, B_\mu}$ is called the \emph{crystal commutor}.
\end{defin}

Note that it is enough to define the crystal commutor
$\sigma_{B_1,B_2} : B_1 \otimes B_2 \to B_2 \otimes B_1$ when $B_1$
and $B_2$ are highest weight crystals since all objects in the
category of $\g$-crystals are unions of these by definition. It was
shown in \cite{HK06,KT06} that for $\g$ a finite-dimensional complex
reductive Lie algebra, the commutor satisfies the cactus relations
\eqref{eq:cactus-relation} and thus endows the category of
$\g$-crystals with the structure of a coboundary (or cactus)
category.  One of the goals of the current paper is to show that
this is true for $\g$ an arbitrary symmetrizable Kac-Moody algebra.


\section{Quiver varieties}
\label{sec:qv}

In this section we introduce the quiver varieties of Lusztig and
Nakajima and the tensor product varieties defined by Malkin and
Nakajima.

\subsection{Lusztig quiver varieties}

Let $I$ be the set of vertices of the Dynkin graph of a Kac-Moody
Lie algebra $\g$ with symmetric Cartan matrix and let $H$ be the set
of pairs consisting of an edge together with an orientation of it.
We call the elements of $H$ \emph{arrows}.  Denote the corresponding
quiver by $Q=(I,H)$. For $h \in H$, let $\inc(h)$ (resp. $\out(h)$)
be the incoming or tip (resp. outgoing or tail) vertex of $h$. We
define the involution $\bar{\ }: H \to H$ to be the function which
takes $h \in H$ to the element of $H$ consisting of the same edge
with opposite orientation. An \emph{orientation} of our graph is a
choice of a subset $\Omega \subset H$ such that $\Omega \cup
\bar{\Omega} = H$ and $\Omega \cap \bar{\Omega} = \emptyset$.  A
\emph{directed path} in $Q$ is a sequence $h_k \dots h_2 h_1$ where
$h_i \in H$ for $1 \le i \le k$ and $\out(h_{i+1}) = \inc(h_i)$ for
$1 \le i \le k-1$. The \emph{length} of such a path is $k$.

Let $\mathcal{V}$ be the category of finite-dimensional $I$-graded
vector spaces $V = \bigoplus_{i \in I} V_i$ over $\C$ with morphisms
being linear maps respecting the grading.  Then $V \in \mathcal{V}$
shall denote that $V$ is an object of $\mathcal{V}$.

Given $V^1,V^2 \in \mathcal{V}$, let
\begin{align*}
E(V^1,V^2) &= \bigoplus_{h \in H} \Hom (V^1_{\out(h)},
V^2_{\inc(h)}) \\
L(V^1,V^2) &= \bigoplus_{i \in I} \Hom (V^1_i,V^2_i).
\end{align*}
For $x = (x_h) \in E(V^1,V^2)$ and $y = (y_h) \in E(V^2,V^3)$,
define $yx \in L(V^1,V^3)$ to be the element with $i$th component
\[
\sum_{h \in H,\, \inc(h)=i} y_h x_{\bar h}
\]
The products $ts, ys, tx$ for $s \in L(V^1,V^2)$ and $t \in
L(V^2,V^3)$ are defined in the obvious way.  For $s \in L(V,V)$, we
define $\tr a = \sum_{i \in I} \tr(a_i)$.

The algebraic group $G_V = \prod_i GL(V_i)$ acts on $E(V,V)$ by $g
\cdot x = g x g^{-1}$. The Lie algebra of $G_V$ is $\mathbf{gl}_{V}
= \prod_i \End(V_i)$.

Define the function $\epsilon : H \to \{-1,1\}$ by
\[
\epsilon(h) = \begin{cases} 1 & \text{if } h \in \Omega \\
-1 & \text{if } h \in {\bar{\Omega}} \end{cases}.
\]
For $x \in E(V^1,V^2)$, define $\epsilon x \in E(V^1,V^2)$ by
$(\epsilon x)_h = \epsilon(h) x_h$ for $h \in H$. Then let
$\left<\cdot,\cdot\right>$ be the nondegenerate, $G_V$-invariant,
symplectic form on $E(V,V)$ with values in $\C$ defined by
\[
\left<x,y\right> = \tr ( (\epsilon x)y).
\]

The moment map associated to the $G_{V}$-action on the symplectic
vector space $E(V,V)$ is the map $\psi : E(V,V) \to \mathbf{gl}_V$
given by
\[
\psi(x) = (\epsilon x)x.
\]
Here we have identified $\mathbf{gl}_V$ with its dual via the trace.

\begin{defin} \label{def:nilpotent}
An element $x \in E(V,V)$ is said to be \emph{nilpotent} if there
exists an $N \ge 1$ such that for any directed path $h_N \dots h_2
h_1$ of length $N$, the composition $x_{h_N} \cdots x_{h_2} x_{h_1}
: V_{\out(h_1)} \to V_{\inc(h_N)}$ is zero.
\end{defin}

Let $\Lambda(V)$ be the set of all nilpotent elements $x \in E(V,V)$
such that $\psi(x)=0$.  Since, up to isomorphism, it depends only on
the graded dimension $\v$ of $V$, we will sometimes denote it
$\Lambda(\v)$.  The variety $\Lambda(V)$ (or $\Lambda(\v)$) is
called a \emph{Lusztig quiver variety}. It was first defined in
\cite{L91}.


\subsection{Nakajima quiver varieties}

For $V, W \in \mathcal{V}$ define
\[
\mathbf{M}(V,W) = E(V,V) \oplus L(W,V) \oplus L(V,W).
\]
The three components of an element of $\mathbf{M}(V,W)$ will
typically be denoted by $x$, $s$, and $t$.  For an $I$-graded
subspace $S \subseteq V$ and $x \in E(V,V)$, we say that $S$ is
\emph{$x$-invariant} if $x_h(S_{\out(h)}) \subseteq S_{\inc(h)}$ for
all $h \in H$.

The group $G_V$ acts on $\mathbf{M}(V,W)$ by
\[
g \cdot (x,s,t) = (gxg^{-1},gs,tg^{-1}).
\]
We have a nondegenerate, $G_V$-invariant, symplectic form on
$\mathbf{M}(V,W)$ defined by
\[
\omega((x,s,t),(x',s',t')) = \tr((\epsilon x)x') + \tr(st' - s't).
\]
The corresponding moment map is given by
\[
\mu(x,s,t) = (\epsilon x)x + st.
\]
Consider the zero set $\mu^{-1}(0)$ of $\mu$.  When we wish to
specify $V$ and $W$, we write $\mu_{V,W}^{-1}(0)$.  This is a (not
necessarily irreducible) affine algebraic variety. We say that a
point $(x,s,t) \in \mu^{-1}(0)$ is \emph{stable} if the only
$I$-graded $x$-invariant subspace of $V$ contained in the kernel of
$t$ is zero.  We denote the set of stable points by $\mu^{-1}(0)^s$.
The action of $G_V$ on $\mu^{-1}(0)^s$ is free and the quotient
\[
\mathfrak{M}(\v,\w) = \mu^{-1}(0)^s/G_V
\]
is a nonsingular quasi-projective variety with symplectic form
induced by $\left<\cdot ,\cdot \right>$.  It is labeled by the
graded dimensions $\v = \bdim V = (\dim V_i)_{i \in I}$ and $\w =
\bdim W = (\dim W_i)_{i \in I}$ of $V$ and $W$ since, up to
isomorphism, it depends only on these dimensions. A $G_V$-orbit
through $(x,s,t)$, considered as a point of $\mathfrak{M}(\v,\w)$,
will be denoted $[x,s,t]$.  We call $\mathfrak{M}(\v,\w)$ a
\emph{Nakajima quiver variety}.  It was originally defined in
\cite{Nak94,Nak98}.

Let $\mathfrak{M}_0(\v,\w) = \mu^{-1}(0) /\!\!/ G_V$ be the affine
algebro-geometric quotient.  That is, it is the affine algebraic
variety whose coordinate ring is the $G_V$-invariant polynomials on
$\mu^{-1}(0)$.  As a set, it consists of the closed $G_V$-orbits in
$\mu^{-1}(0)$.  We have a projective morphism $\pi :
\mathfrak{M}(\v,\w) \to \mathfrak{M}_0(\v,\w)$ which sends $[x,s,t]$
to the unique closed orbit contained in the closure of the orbit
$G_V \cdot (x,s,t)$.  We then define
\[
\mathfrak{L}(\v,\w) = \pi^{-1}(0).
\]
It is a lagrangian subvariety of $\mathfrak{M}(\v,\w)$.  Let
\[
\mathfrak{M}(\w) = \bigsqcup_\v \mathfrak{M}(\v,\w),\quad
\mathfrak{L}(\w) = \bigsqcup_\v \mathfrak{L}(\v,\w).
\]


\subsection{Tensor product quiver varieties}

Let $W^i \in \mathcal{V}$, $i=1,2,\dots,n$, with graded dimensions
$\w^i$ and $V \in \mathcal{V}$ with graded dimension $\v$.  Set $W =
\bigoplus_{i=1}^n W^i$ and $W^{j,k} = \bigoplus_{i=j}^k W^i$ for $1
\le j \le k \le n$. We adopt the convention that $W^{j,k}=0$ if $j
> k$.  The group $G_W$ acts on $\mathbf{M}(V,W)$ by
\[
g * (x,s,t) = (x,sg^{-1},gt).
\]
This commutes with the action of $G_V$ and thus induces an action of
$G_W$ on $\mathfrak{M}(\v,\w)$ and $\mathfrak{M}_0(\v,\w)$ and the
map $\pi$ is $G_W$-equivariant. Define a one-parameter subgroup
$\lambda : \C^* \to G_W$ by
\[
\lambda(z) = \id_{W^1} \oplus z\, \id_{W^2} \oplus \cdots \oplus
z^{n-1}\, \id_{W^n} \in \prod_{i=1}^n GL(W^i) \subseteq G_W.
\]
Let $\mathfrak{M}(\v,\w)^{\lambda(\C^*)}$ denote the fixed point set
of $\mathfrak{M}(\v,\w)$.

\begin{lem}[{\cite[Lemma~3.2]{Nak01}}]
We have
\[
\mathfrak{M}(\v,\w)^{\lambda(\C^*)} \cong \bigsqcup_{\sum \v^i = \v}
\left( \prod_{i=1}^n \mathfrak{M}(\v^i,\w^i) \right)
\]
where the union is over all ordered $n$-tuples $(\v^1,
\v^2, \dots, \v^n)$ such that $\sum \v^i = \v$.
\end{lem}
Taking the union over all possible $\v$ yields
\[
\mathfrak{M}(\w)^{\lambda(\C^*)} \cong \prod_{i=1}^n
\mathfrak{M}(\w^i).
\]
Define the \emph{tensor product quiver variety}
\begin{gather*}
\T = \left\{ [x,s,t] \in \mathfrak{M}(\w)\ |\ \lim_{z \to 0}
\lambda(z) * [x,s,t] \in \mathfrak{L}(\w^1) \times \dots \times
\mathfrak{L}(\w^n) \right\}, \text{ and}\\
\vT = \T \cap \mathfrak{M}(\v,\w).
\end{gather*}
Note that the limit in the above definition does not always exist.
As shown in \cite[Lemma~3.6]{Nak01}, $\T$ is a closed subvariety of
$\mathfrak{M}(\w)$ and $\vT$ is a closed subvariety of
$\mathfrak{M}(\v,\w)$.  By \cite[Lemma~3.5]{Nak01}, we also have
\[
\T = \left\{ [x,s,t] \in \mathfrak{M}(\w)\ |\ \lim_{z \to 0}
\lambda(z) * \pi ([x,s,t])=0 \right\}.
\]

\begin{example}[$\mathfrak{g} = \mathfrak{sl}_2$]
If $\mathfrak{g}=\mathfrak{sl}_2$, the corresponding quiver has one
vertex and no edges.  Therefore $x=0$, the stability condition
forces $t$ to be injective, and the moment map condition implies
$st=0$. Then the map
\[
    (0,s,t) \mapsto (\im t, ts)
\]
identifies $\mathfrak{M}(v,w)$ with
\[
    \{(S,\chi) \in Gr(v,w) \times \End W\ |\ \chi(S) =
    0,\ \chi(W) \subseteq S\} \cong T^*Gr(v,w),
\]
where $Gr(v,w)$ is the Grassmannian of dimension $v$ planes in the
$w$-dimensional space $W$.  If we fix a decomposition $W =
\bigoplus_{i=1}^n W^i$, then $\mathfrak{T}(w^1, \dots, w^n)$ is the
subvariety of $\mathfrak{M}(v,w)$ consisting of those $(S,\chi) \in
Gr(v,w) \times \End W$ for which $\chi(W^i) \subseteq W^{i+1} \oplus
\dots \oplus W_n$.
\end{example}


\section{Geometric realizations of crystals}
\label{sec:geometric-crystals}

In this section we recall the construction of the crystals
$B_\infty$ and $B_{\lambda_1} \otimes \dots \otimes B_{\lambda_n}$
on sets of irreducible components of quiver varieties.  We also
describe a geometric realization of Kashiwara's involution.

\subsection{Geometric realization of $B_\infty$}
\label{sec:geom-real-Binf}

We briefly review here the geometric realization of the crystal
$B_\infty$ defined by Kashiwara and Saito \cite{KS97}.  We identify
$(\Z_{\ge 0})^I$ with the negative root lattice of $\g$ by
identifying $\v = (\v_i)$ with $-\sum_{i \in I} \v_i \alpha_i$ where
$\alpha_i$ are the simple roots of $\g$.  For each $\v \in (\Z_{\ge
0})^I$, choose an $I$-graded vector space $V(\v)$ of graded
dimension $\v$. Then let $\Lambda(\v) = \Lambda(V(\v))$. Let
$\Lambda(\v',\v)$ be the variety of triples $(x,\phi',\bar \phi)$
such that $x \in \Lambda(\v)$ and $\phi'=(\phi'_i)$, $\bar \phi =
(\bar \phi_i)$ give an exact sequence
\[
0 \to V(\v')_i \stackrel{\phi'_i}{\longrightarrow} V(\v)_i
\stackrel{\bar \phi_i}{\longrightarrow} V(\v - \v')_i \to 0
\]
for each $i \in I$ and $\im \phi'$ is $x$-invariant.  Then $x$
induces a map $x' \in \Lambda(\v')$ and so we have the following
diagram
\begin{equation} \label{eq:LQV-two-projs}
\Lambda(\v') \stackrel{q_1}{\longleftarrow} \Lambda(\v',\v)
\stackrel{q_2}{\longrightarrow} \Lambda(\v)
\end{equation}
where $q_1(x,\phi',\bar \phi) = x'$ and $q_2(x,\phi',\bar \phi) =
x$.

For $x \in \Lambda(\v)$ and $i \in I$, let
\[
\varepsilon_i(x) = \dim \Coker \left( \bigoplus_{h\, :\, \inc(h)=i}
V(\v)_{\out(h)} \stackrel{(x_h)}{\longrightarrow} V(\v)_i \right)
\]
and for $i \in I$ and $c \in \Z_{\ge 0}$ let
\[
\Lambda(\v)_{i,c} = \{x \in \Lambda(\v)\ |\ \varepsilon_i(x)=c\}.
\]
Let $B(\v,\infty)$ be the set of irreducible components of
$\Lambda(\v)$ and for $X \in B(\v,\infty)$, define $\varepsilon_i(X)
= \varepsilon_i(x)$ for a generic point $x$ of $X$.  For $i \in I$
and $c \in \Z_{\ge 0}$, let
\[
B(\v,\infty)_{i,c} = \{X \in B(\v,\infty)\ |\ \varepsilon_i(X)=c\}.
\]
Then \eqref{eq:LQV-two-projs} induces an isomorphism (see
\cite[Prop~5.2.4]{KS97})
\begin{equation} \label{eq:LQV-irrcomp-isom}
B(\v+c\alpha_i,\infty)_{i,0} \cong B(\v,\infty)_{i,c}.
\end{equation}

Let $B(\infty) = \bigsqcup_{\v} B(\v,\infty)$.  We define crystal
operators on $B(\infty)$ as follows.  Suppose $\Lambda' \in B(\v +
c\alpha_i,\infty)_{i,0}$ corresponds to $\Lambda \in
B(\v,\infty)_{i,c}$ by the isomorphism \eqref{eq:LQV-irrcomp-isom}.
Define
\begin{gather*}
\kf_i^c : B(\v + c\alpha_i, \infty)_{i,0} \to B(\v,\infty)_{i,c},
\quad \kf_i^c(\Lambda') = \Lambda, \\
\ke_i^c : B(\v, \infty)_{i,c} \to B(\v + c\alpha_i,\infty)_{i,0},
\quad \ke_i^c(\Lambda) = \Lambda'.
\end{gather*}
For $c > 0$ we then define $\ke_i : B(\infty) \to B(\infty) \sqcup
\{0\}$ by
\[
\ke_i : B(\v,\infty)_{i,c} \stackrel{\ke_i^c}{\longrightarrow} B(\v
+ c\alpha_i,\infty)_{i,0} \stackrel{\kf_i^{c-1}}{\longrightarrow}
B(\v + \alpha_i,\infty)_{i,c-1}
\]
and let $\ke_i(X)=0$ for $X \in B(\v,\infty)_{i,0}$.  Also define
\[
\kf_i : B(\infty) \to B(\infty), \quad \kf_i : B(\v,\infty)_{i,c}
\stackrel{\ke_i^c}{\longrightarrow} B(\v+c\alpha_i,\infty)_{i,0}
\stackrel{\kf_i^{c+1}}{\longrightarrow}
B(\v-\alpha_i,\infty)_{i,c+1}.
\]
Then $\ke_i^c$ and $\kf_i^c$ can be considered the $c$th powers of
$\ke_i$ and $\kf_i$ respectively.  If $P$ is the weight lattice of
$\g$, we also define
\begin{gather*}
\wt : B(\infty) \to P; \quad \wt(X) = \v \text{ for } X \in
B(\v,\infty),\\
\varphi_i(X) = \varepsilon_i(X) + \left<h_i,\wt(X)\right>.
\end{gather*}

\begin{prop}[{\cite[Thm~5.2.6, Thm~5.3.2]{KS97}}]
The above definitions endow $B(\infty)$ with the structure of a
$\g$-crystal and $B(\infty) \cong B_\infty$ as $\g$-crystals.
\end{prop}

For an element $b \in B_\infty$, let $X_b$ denote the corresponding
element of $B(\infty)$.


\subsection{Geometric realization of Kashiwara's involution}
\label{sec:geom-kash-involution}

We recall here a geometric realization, introduced by Kashiwara and
Saito \cite{KS97}, of the involution described in
Section~\ref{sec:kash-involution}. Fix a nondegenerate Hermitian
form on $V(\v)_i$ for all $i$ and $\v$.  Then $x \mapsto x^\dag$,
where $\dag$ denotes the Hermitian adjoint, gives an automorphism of
$E(V(\v),V(\v))$ and $\Lambda(\v)$ is invariant under this
automorphism.  This induces an involution of $B(\v,\infty)$ which we
denote by $*$.  Since $\Lambda(\v)$ is $G_{V(\v)}$-invariant, the
involution $*$ does not depend on our choice of Hermitian forms.  It
was shown in \cite{KS97} that $*$ corresponds to Kashiwara's
involution under the isomorphism $B(\infty) \cong B_\infty$.  That
is, $X_b^* = X_{b^*}$ for all $b \in B_\infty$.  Note that in
\cite{KS97}, an isomorphism between $V(\v)_i$ and its dual was
chosen for each $\v$ and $i$ and the transpose, rather than the
Hermitian adjoint, was used to realize Kashiwara's involution.
Fixing a real form of each $V(\v)_i$ (i.e. a real vector space
$V(\v)_i^{\mathbb{R}}$ such that $V(\v)_i = V(\v)_i^{\mathbb{R}}
\otimes_{\mathbb{R}} \C$), our Hermitian form yields a nondegenerate
bilinear form given by $(u,v) \mapsto \left<u, \kappa(v)\right>$
where $\left< \cdot, \cdot \right>$ is the Hermitian form and
$\kappa : w \otimes z \mapsto w \otimes \bar z$, $w \in
V(\v)_i^{\mathbb{R}}$, $z \in \C$, is the conjugation determined by
the real form. This gives an identification of $V(\v)_i$ with its
dual for each $i$ and Hermitian adjoint corresponds to transpose.
In what follows, we will often write $V$ for $V(\v)$, when it will
cause no confusion, to simplify notation.


\subsection{Geometric realization of $B_{\lambda_1} \otimes \dots
\otimes B_{\lambda_n}$} \label{sec:TPQV-crystal}

Malkin \cite{Mal03} and Nakajima \cite{Nak01} have endowed the set
of irreducible components of the tensor product quiver variety with
the structure of a $\g$-crystal.  We briefly recall the construction
here.

Let
\begin{multline}
\mathfrak{T}(\v',\v ; \w^1,\dots, \w^n) = \{(x,s,t,S)\ |\ [x,s,t]
\in \vT,\,\\ \im s \subseteq S \subseteq V,\, S \text{ is
$x$-invariant},\, \bdim S = \v' \}/G_V.
\end{multline}
Then we have the diagram
\begin{equation} \label{eq:TPQV-two-projs}
\mathfrak{T}(\v' ; \w^1,\dots,\w^n) \stackrel{q'_1}{\longleftarrow}
\mathfrak{T}(\v',\v ; \w^1,\dots,\w^n)
\stackrel{q'_2}{\longrightarrow} \vT
\end{equation}
where $q'_1(G_V \cdot (x,s,t,S)) = [x^S,s^{W,S},t^{S,W}]$ and
$q'_2(G_V \cdot (x,s,t,S)) = [x,s,t]$.  Here $x^S$ and $t^{S,W}$
denote the restriction of $x$ and $t$ to $S$ respectively and
$s^{W,S}$ is the map $s$ viewed as a map into $S$.

For $[x,s,t] \in \vT$ and $i \in I$, define
\[
\varepsilon_i([x,s,t]) = \dim (V_i/\im \tau_i),
\]
where
\[
\tau_i = \tau_{i,(x,s,t)} = \sum_{h\, :\, \inc(h)=i} \epsilon(h) x_h
+ s_i : \bigoplus_{h\, :\, \inc(h)=i} V_{\out(h)} \oplus W_i \to
V_i.
\]
For $c \in \Z_{\ge 0}$ let
\[
\vT_{i,c} = \{[x,s,t] \in \vT\ |\ \varepsilon_i([x,s,t]) = c\}.
\]
Then $\vT_{i,c}$ is a locally closed subvariety of $\vT$.  Let
$B(\v;\w^1, \dots, \w^n)$ denote the set of irreducible components
of $\vT$ and let $B(\w^1, \dots, \w^n) = \bigsqcup_{\v} B(\v;\w^1,
\dots, \w^n)$. For $X \in B(\v;\w^1, \dots, \w^n)$ define
$\varepsilon_i(X) = \varepsilon_i([x,s,t])$ for a generic point
$[x,s,t]$ of $X$.  For $i \in I$ and $c \in \Z_{\ge 0}$, let
\[
B(\v;\w^1, \dots, \w^n)_{i,c} = \{X \in B(\v;\w^1, \dots, \w^n)\ |\
\varepsilon_i(X)=c\}.
\]
For $\w \in (\Z_{\ge 0})^I$, let $\lambda_\w = \sum_i \w_i \omega_i$
where the $\omega_i$ are the fundamental weights of $\g$.  Then
define
\begin{gather*}
\wt : B(\w^1, \dots, \w^n) \to P,\quad \wt(X) = \lambda_\w + \v
\text{ for } X \in B(\v;\w^1, \dots, \w^n), \\
\varphi_i(X) = \varepsilon_i(X) + \left< h_i, \wt(X) \right>.
\end{gather*}
Note that for $X \in B(\v;\w^1, \dots, \w^n)_{i,c}$,
\[
\varphi_i(X) = c + \left< h_i, \lambda_\w + \v \right> = \dim (\ker
\tau_i / \im \gamma_i), \quad \text{where}\quad  \gamma_i =
\bigoplus_{h\, :\, \inc(h)=i} x_{\bar h} \oplus t_i.
\]

The maps \eqref{eq:TPQV-two-projs} induce an isomorphism (see
\cite[{\S}4]{Nak01})
\begin{equation} \label{eq:TPQV-irrcomp-isom}
B(\v + c\alpha_i;\w^1,\dots,\w^n)_{i,0} \cong B(\v;\w^1, \dots,
\w^n)_{i,c}.
\end{equation}
We define crystal operators on $B(\w^1, \dots, \w^n)$ as follows.
Let $X' \in B(\v+c\alpha_i;\w^1, \dots, \w^n)_{i,0}$ correspond to
$X \in B(\v;\w^1, \dots, \w^n)_{i,c}$ under the isomorphism
\eqref{eq:TPQV-irrcomp-isom}. Then define
\begin{gather*}
\kf_i^c : B(\v + c\alpha_i;\w^1, \dots, \w^n)_{i,0} \to B(\v;\w^1,
\dots, \w^n)_{i,c}, \quad \kf_i^c(X')=X, \\
\ke_i^c : B(\v;\w^1, \dots, \w^n)_{i,c} \to B(\v + c\alpha_i;\w^1,
\dots, \w^n)_{i,0},\quad \ke_i^c(X)=X'.
\end{gather*}
For $c > 0$, we then define $\ke_i: B(\w^1, \dots, \w^n) \to B(\w^1,
\dots, \w^n)$ by
\[
\ke_i : B(\v;\w^1, \dots, \w^n)_{i,c}
\stackrel{\ke_i^c}{\longrightarrow} B(\v + c\alpha_i; \w^1, \dots,
\w^n)_{i,0} \stackrel{\kf_i^{c-1}}{\longrightarrow} B(\v +
\alpha_i;\w^1, \dots, \w^n)_{i,c-1},
\]
and set $\ke_i(X) = 0$ for $X \in B(\v;\w^1, \dots, \w^n)_{i,0}$.
For $c > -\left< h_i, \lambda_\w + \v \right>$, let
\[
\kf_i : B(\v;\w^1, \dots, \w^n)_{i,c}
\stackrel{\ke_i^c}{\longrightarrow} B(\v + c\alpha_i; \w^1, \dots,
\w^n)_{i,0} \stackrel{\kf_i^{c+1}}{\longrightarrow} B(\v-\alpha_i;
\w^1, \dots, \w^n)_{i,c+1},
\]
and set $\kf_i(X)=0$ for $X \in B(\v;\w^1, \dots, \w^n)_{i,c}$ with
$c \le - \left< h_i, \lambda_\w + \v \right>$.  The maps $\ke_i^c$
and $\kf_i^c$ defined above can be considered the $c$th powers of
$\ke_i$ and $\kf_i$ respectively.

\begin{prop}[{\cite[Prop~4.3, Thm~4.6, Cor~4.7, {\S}7]{Nak01}}]
\label{prop:TPQV-crystal}
The above definitions endow $B(\w^1, \dots, \w^n)$ with the
structure of a $\g$-crystal and $B(\w^1, \dots, \w^n) \cong
B_{\lambda_{\w^1}} \otimes \dots \otimes B_{\lambda_{\w^n}}$ as
$\g$-crystals.
\end{prop}

For $b \in B_{\lambda_{\w^1}} \otimes \dots \otimes
B_{\lambda_{\w^n}}$, let $Y_b \in B(\w^1, \dots, \w^n)$ denote the
corresponding irreducible component of the tensor product quiver
variety $\T$.


\subsection{Fiber bundles and crystal isomorphisms}

Let $\tT$ and $\tvT$ denote the inverse images of $\T$ and $\vT$
(respectively) under the natural projection $\mu^{-1}(0)^s \to
\mathfrak{M}(\w)$. That is,
\begin{align*}
\tT &= \{(x,s,t) \in \mu^{-1}(0)^s\ |\ [x,s,t] \in \T\}, \\
\tvT &= \{(x,s,t) \in \mu^{-1}(0)^s\ |\ [x,s,t] \in \vT\}.
\end{align*}
Since the aforementioned projection is a principle $G_V$-bundle, the
irreducible components of $\tT$ are in natural one-to-one
correspondence with the irreducible components of $\T$.  Note that
the $n=1$ case reduces to $\mathfrak{T}(\w) = \mathfrak{L}(\w)$ and
we define $\widetilde{\mathfrak{L}}(\w) =
\widetilde{\mathfrak{T}}(\w)$. Let $\tilde Y_b$ denote the
irreducible component of $\tT$ corresponding to the irreducible
component $Y_b$ of $\T$.

It will be useful to have a slightly more concrete description of
$\T$. It is shown in \cite[Prop~3.8]{Nak01} (while only the case
$n=2$ is considered there, the generalization to higher $n$ is
straightforward) that $\T$ decomposes as a disjoint union
\[
\T = \bigsqcup_{\v^1,\v^2,\dots,\v^n} \mathfrak{T}(\v^1, \dots, \v^n
| \w^1, \dots, \w^n)
\]
where
\begin{multline*}
\mathfrak{T}(\v^1, \dots, \v^n
| \w^1, \dots, \w^n) \\
\stackrel{\text{def}}{=} \left\{ [x,s,t] \left| \lim_{z \to 0}
\lambda(z) * [x,s,t] \in \mathfrak{L}(\v^1,\w^1) \times
\mathfrak{L}(\v^2,\w^2) \times \dots \times \mathfrak{L}(\v^n,\w^n)
\right. \right\}.
\end{multline*}
These are the Bialynicki-Birula decompositions of $\T$. The map
\begin{multline} \label{eq:fiber-bundle}
\mathfrak{T}(\v^1, \dots, \v^n | \w^1, \dots, \w^n) \ni [x,s,t]
\mapsto \lim_{z \to 0} \lambda(z) * [x,s,t] \\
\in \mathfrak{L}(\v^1,\w^1) \times \mathfrak{L}(\v^2,\w^2) \times
\dots \times \mathfrak{L}(\v^n,\w^n)
\end{multline}
is a fiber bundle with affine fibers.  By a generalization of the
results of \cite[Prop~3.8,~Prop~3.15]{Nak01} to more than two
factors, these fiber bundles identify the irreducible components of
$\T$ with the irreducible components of $\mathfrak{L}(\w^1) \times
\dots \times \mathfrak{L}(\w^n)$ and this identification is an
isomorphism of crystals \cite[Thm~4.6]{Nak01}. Here we use the
tensor product rule and the crystal structure on each
$\mathfrak{L}(\w^i)$ to give a crystal structure to
$\mathfrak{L}(\w^1) \times \dots \times \mathfrak{L}(\w^n)$. In
general, the crystal $B_{\lambda_{\w^1}} \otimes \dots \otimes
B_{\lambda_{\w^n}}$ has nontrivial automorphisms.  Therefore, the
isomorphism of Proposition~\ref{prop:TPQV-crystal} is not
necessarily unique. However, each $B_{\lambda_{\w^i}}$ has no
nontrivial automorphisms since it is generated by a single highest
weight element. Therefore we can use the identification of $B(\w^1,
\dots, \w^n)$ with $B(\w^1) \times \dots \times B(\w^n)$ induced by
\eqref{eq:fiber-bundle} and the unique isomorphisms $B(\w^i) \cong
B_{\lambda_{\w^i}}$ to fix an isomorphism
\begin{equation} \label{eq:phi-isom}
\phi :  B(\w^1, \dots, \w^n) \stackrel{\cong}{\longrightarrow}
B_{\lambda_{\w^1}} \otimes \dots \otimes B_{\lambda_{\w^n}}.
\end{equation}

\begin{lem} \label{lem:one-param-limits}
Suppose $[x,s,t] \in \mathfrak{T}(\v^1, \dots, \v^n | \w^1, \dots,
\w^n)$.  That is,
\[
\lim_{z \to 0} \lambda(z) * [x,s,t] \in \mathfrak{L}(\v^1,\w^1)
\times \mathfrak{L}(\v^2,\w^2) \times \dots \times
\mathfrak{L}(\v^n;\w^n).
\]
Then there exists a $d \in \Z_{>0}$, representatives $(x^j, s^j,
t^j) \in \tilde {\mathfrak{L}}(\v^j, \w^j)$ and a one-parameter
subgroup $\rho : \C^* \to G_V$ such that
\begin{equation} \label{eq:one-param-limit-quotient}
\lim_{z \to 0} \lambda(z)
* [x,s,t] = ([x^1,s^1,t^1], \dots, [x^n,s^n,t^n])
\end{equation}
and
\begin{equation} \label{eq:one-param-limit-upstairs}
\lim_{z \to 0} \lambda(z^d) * \rho(z) \cdot (x,s,t) = (x^1 \oplus
\dots \oplus x^n, s^1 \oplus \dots \oplus s^n, t^1 \oplus \dots
\oplus t^n).
\end{equation}
\end{lem}

\begin{proof}
Let
\[
    [x',s',t'] = \lim_{z \to 0} \lambda(z) * [x,s,t]
\]
and fix a representative $(x',s',t') \in \mu^{-1}(0)^s \subseteq
\mathbf{M}(V,W)$.  Then we can write
\[
(x', s', t') = (x^1 \oplus \dots \oplus x^n, s^1 \oplus \dots \oplus
s^n, t^1 \oplus \dots \oplus t^n)
\]
for some $(x^j, s^j, t^j) \in \tilde {\mathfrak{L}}(\v^j, \w^j)$, $1
\le j \le n$.  By the geometric invariant theory (see, for example,
\cite[Definition~1.7]{MF94}), there exists an affine $G_V$-invariant
neighborhood $U$ of $(x',s',t')$ in $\mathbf{M}(V,W)$ such that the
orbit $G_V \cdot (x',s',t')$ is closed in $U$.  We may assume that
$(x,s,t)$ is in $U$.  Consider the action of the reductive group
$\C^* \times G_V$ on $\mathbf{M}(V,W)$ given by
\[
((z,g), (x,s,t)) \mapsto (z,g) \star (x,s,t)
\stackrel{\text{def}}{=} \lambda(z)
* g \cdot (x,s,t).
\]
By hypothesis, $G_V \cdot (x',s',t')$  meets the closure of the
orbit $(\C^* \times G_V) \star (x,s,t)$.  Therefore, by a version of
the Hilbert criterion (see \cite[Thm~1.4]{Kem78}), since $U$ is
affine, there exists a one-parameter subgroup $(\rho', \rho) : \C^*
\to \C^* \times G_V$ such that $\lim_{z \to 0} (\rho', \rho)(z)
\star (x,s,t)$ exists and is contained in $G_V \cdot (x',s',t')$. By
modification of the one-parameter subgroup (or representative
$(x',s',t')$), we may assume that
\[
\lim_{z \to 0} (\rho', \rho)(z) \star (x,s,t) = (x',s',t').
\]
We may also assume that $\rho' (z) \to 0$ as $z \to 0$. That is
$\rho' (z) = z^d$ for some $d \in \Z_{>0}$.  Then we have
\[
\lim_{z \to 0} \lambda(z^d) * \rho(z) \cdot (x,s,t) = \lim_{z \to 0}
(\rho', \rho)(z) \star (x,s,t) = (x', s', t')
\]
as desired.
\end{proof}

Note that if we have a flag of $I$-graded spaces
\[
0 = V^{n+1} \subseteq V^n \subseteq V^{n-1} \subseteq \dots
\subseteq V^1 = V,
\]
with
\[
x (V^i) \subseteq V^i,\quad s(W^i) \subseteq V^i,\quad t(V^i)
\subseteq W^{i,n},\quad 1 \le i \le n,
\]
then $x$, $s$ and $t$ induce maps
\[
x^{V^i/V^{i+1}} : V^i/V^{i+1} \to V^i/V^{i+1},\quad
s^{W^i,V^i/V^{i+1}} : W^i \to V^i/V^{i+1},\quad t^{V^i/V^{i+1},W^i}
: V^i/V^{i+1} \to W^i
\]

\begin{prop} \label{prop:fiber-map}
The set $\mathfrak{T}(\v^1, \dots, \v^n | \w^1, \dots, \w^n)$
consists of those $[x,s,t]$ in $\T$ such that there exists a flag of
$I$-graded spaces
\[
0 = V^{n+1} \subseteq V^n \subseteq V^{n-1} \subseteq \dots
\subseteq V^1 = V,\quad \bdim V^i/V^{i+1} = \v_i,
\]
with
\[
x (V^i) \subseteq V^i,\quad s(W^i) \subseteq V^i,\quad t(V^i)
\subseteq W^{i,n},
\]
and
\[
(x^{V^i/V^{i+1}}, s^{W^i, V^i/V^{i+1}}, t^{V^i/V^{i+1},W^i}) \in
\widetilde{\mathfrak{L}}(\v^i, \w^i),\quad 1 \le i \le n.
\]
\end{prop}

\begin{proof}
By Lemma~\ref{lem:one-param-limits}, there exists a positive integer
$d$, representatives $(x^j, s^j, t^j) \in \tilde
{\mathfrak{L}}(\v^j, \w^j)$, and a one-parameter subgroup $\rho :
\C^* \to G_V$ such that \eqref{eq:one-param-limit-quotient} and
\eqref{eq:one-param-limit-upstairs} hold.  Denote the $\rho$-weight
space decomposition of $V_i$, $i \in I$, by
\[
V_i = \bigoplus_{m \in \Z} V_i^{(m)},\quad V_i^{(m)} = \{v \in V_i\
|\ \rho(z)_i(v) = z^m v, z \in \C^*\}.
\]
The sum $V^{(m)} = \bigoplus_{i \in I} V_i^{(m)}$ is an $I$-graded
subspace of $V$.  Then \eqref{eq:one-param-limit-upstairs} implies
\[
x(V^{(k)}) \subseteq \bigoplus_{m \ge k} V^{(m)},\quad s(W^k)
\subseteq \bigoplus_{m \ge d(k-1)} V^{(m)},\quad
t\left(V^{(d(k-1))}\right) \subseteq W^{k,n}.
\]
The stability condition then implies that $V^{(l)}=0$ for $l \ge dn$
and the flag given by
\[
V^k = \bigoplus_{m \ge d(k-1)} V^{(m)},\quad 1 \le k \le n.
\]
satisfies the conditions of the proposition.

Conversely, suppose that for some $[x,s,t] \in \T$, a flag with the
given properties exists. For each $i \in I$, choose a decomposition
\begin{equation} \label{eq:Vi-decomp}
V_i = \bigoplus_{m=1}^n V_i^{(m)}
\end{equation}
such that
\[
V^k_i = \bigoplus_{m \ge k} V^{(m)}_i,\quad 1 \le k \le n.
\]
Then define a one-parameter subgroup $\rho : \C^* \to G_V$ by
\[
\rho(z)_i = \sum_{m=1}^n z^{m-1} \id_{V_i^{(m)}}.
\]
Then it is easily seen that
\[
\lim_{z \to 0} \lambda(z) * \rho(z) \cdot (x,s,t) = (x^{V^{(1)}}
\oplus \dots \oplus x^{V^{(n)}}, s^{W^1, V^{(1)}} \oplus \dots
\oplus s^{W^n, V^{(n)}}, t^{V^{(1)}, W^1} \oplus \dots \oplus
t^{V^{(n)}, W^n})
\]
where $x^{V^{(k)}}$ denotes the restriction of $x$ to $V^{(k)}$
composed with the projection to $V^{(k)}$ (according to the
decomposition given in \eqref{eq:Vi-decomp}). The maps $s^{W^k,
V^k}$ and $t^{V^k, W^k}$ are defined similarly. Thus
\[
\lim_{z \to 0} \lambda(z) * [x,s,t] = ([x^1,s^1,t^1], \dots,
[x^n,s^n,t^n]).
\]
where
\[
[x^k, s^k, t^k] = [x^{V^k/V^{k+1}}, s^{W^k, V^k/V^{k+1}}, t^{V^k/
V^{k+1}, W^k}] \in \mathfrak{L}(\v^k,\w^k).
\]
\end{proof}

We define
\[
\widetilde{\mathfrak{T}}(\v^1,\dots,\v^n | \w^1, \dots, \w^n) =
\left\{ (x,s,t)\ |\ [x,s,t] \in \mathfrak{T}(\v^1, \dots, \v^n |
\w^1, \dots, \w^n)\right\}.
\]

For $0 = p_0 < p_1 \le p_2 \le \dots \le p_k \le n$, define a more
general one-parameter subgroup $\lambda^{(p_1,\dots,p_k)} : \C^* \to
G_W$ by
\[
\lambda^{(p_1,\dots,p_k)}(z) = \id_{W^{1,p_1}} \oplus z \id_{W^{p_1
+ 1,p_2}} \oplus \dots \oplus z^k \id_{W^{p_k+1,n}}.
\]
Then, as above, $\T$ decomposes as a disjoint union
\[
\T = \bigsqcup_{\v^1,\dots,\v^{k+1}} \mathfrak{T}^{(p_1, \dots,
p_k)}(\v^1, \dots, \v^{k+1} | \w^1, \dots, \w^n)
\]
where
\begin{multline}
\mathfrak{T}^{(p_1, \dots, p_k)}(\v^1, \dots, \v^{k+1} | \w^1,
\dots, \w^n) \stackrel{\text{def}}{=} \left\{ [x,s,t] \left| \lim_{z
\to 0} \lambda^{(p_1, \dots, p_k)}(z) * [x,s,t] \in \right. \right.
\\
\left. \mathfrak{T}(\v^1 ; \w^1, \dots, \w^{p_1}) \times \dots
\times \mathfrak{T}(\v^{k+1} ; \w^{p_k+1}, \dots, \w^n) \right\}.
\end{multline}
The map
\begin{multline} \label{eq:p1k-fiber-bundle}
\mathfrak{T}^{(p_1, \dots, p_k)}(\v^1, \dots, \v^{k+1} | \w^1,
\dots, \w^n) \ni [x,s,t]
\mapsto \lim_{z \to 0} \lambda^{(p_1, \dots, p_k)}(z) * [x,s,t] \in \\
\mathfrak{T}(\v^1 ; \w^1, \dots, \w^{p_1}) \times \dots \times
\mathfrak{T}(\v^{k+1} ; \w^{p_k+1}, \dots, \w^n)
\end{multline}
is a fiber bundle with affine fibers.  These fiber bundles identify
the irreducible components of $\T$ with the irreducible components
of $\mathfrak{T}(\w^1, \dots, \w^{p_1}) \times \dots \times
\mathfrak{T}(\w^{p_k+1}, \dots, \w^n)$ and this identification is an
isomorphism of crystals.  We then have the following generalization
of Proposition~\ref{prop:fiber-map}.

\begin{prop} \label{prop:fiber-map-general}
The set $\mathfrak{T}^{(p_1, \dots, p_k)}(\v^1, \dots, \v^{k+1} |
\w^1, \dots, \w^n)$ consists of those $[x,s,t]$ in $\T$ such that
there exists a flag of $I$-graded spaces
\[
0 = V^{k+2} \subseteq V^{k+1} \subseteq V^k \subseteq \dots
\subseteq V^1 = V,\quad \bdim V^i/V^{i+1} = \v^i,
\]
with
\[
x (V^i) \subseteq V^i,\quad s(W^{p_{i-1}+1,p_i}) \subseteq V^i,\quad
t(V^i) \subseteq W^{p_{i-1}+1,n},
\]
and
\[
(x^{V^i/V^{i+1}}, s^{W^{p_{i-1}+1,p_i}, V^i/V^{i+1}},
t^{V^i/V^{i+1},W^{p_{i-1}+1,p_i}}) \in {\widetilde
{\mathfrak{T}}}(\v^i ; \w^{p_{i-1}+1}, \dots, \w^{p_i}),\quad 1 \le
i \le n.
\]
\end{prop}
\begin{proof}
This follows from Proposition~\ref{prop:fiber-map}.  The details are
left to the reader.
\end{proof}

\begin{prop}
\label{prop:x-nilpotent} If $(x,s,t) \in
\widetilde{\mathfrak{T}}(\w^1;\dots;\w^n)$ (equivalently, $[x,s,t]
\in \mathfrak{T}(\w^1; \dots; \w^n)$) then $x$ is nilpotent.
\end{prop}
\begin{proof}
This follows from Proposition~\ref{prop:fiber-map} and the fact that
$x$ is nilpotent for all $(x,s,t) \in
\widetilde{\mathfrak{L}}(\v^i,\w^i)$ (see the proof of
\cite[Lemma~5.9]{Nak94}).
\end{proof}


\section{A geometric commutor}
\label{sec:geom-commutor}

In this section, we give precise characterizations of the
irreducible components of the tensor product quiver variety.  We
then examine the action of the crystal commutor in terms of these
characterizations.  This will enable us to prove that the commutor
satisfies the cactus relation in Section~\ref{sec:cactus-relation}.

\subsection{Characterization of irreducible components}

For an arbitrary element $b \in B_{\lambda_1} \otimes \dots \otimes
B_{\lambda_n}$, let $\hw b$ be the unique highest weight element in
the connected component of the crystal graph of $B_{\lambda_1}
\otimes \dots \otimes B_{\lambda_n}$ containing $b$.  If $b'$ is a
highest weight element in $B_{\lambda_1} \otimes \dots \otimes
B_{\lambda_n}$, then we can identify the connected component
containing $b'$ with some $B_\lambda$ and $b'$ corresponds to
$b_\lambda$ under this identification.  For $b \in B_\infty$ and
$b'$ as above, we then define
\[
    \tilde b b' = \begin{cases} b'' & \text{if }
    \iota^\infty_\lambda(b'') = b, \\
    0 & \text{if } \not \exists \ b'' \text{ such that }
    \iota^\infty_\lambda(b'') = b, \end{cases}
\]
(this is well defined since $\iota^\infty_\lambda$ is injective).
Equivalently, $\tilde b b' = \Phi_\lambda(b)$, where $\Phi_\lambda :
B_\infty \to B_\lambda \sqcup \{0\}$ is the natural crystal morphism
sending $b_\infty$ to $b_\lambda$.  We view $\tilde b b'$ as an
element of $B_{\lambda_1} \otimes \dots \otimes B_{\lambda_n}$ via
the above identification. Note that for $b \in B_\lambda$, $b_1 =
\iota^\infty(b)$, we have $b = \tilde b_1 b_\lambda$.

\begin{defin} \label{def:bp}
\begin{enumerate}
\item For $k \ge 1$ and
$0 = p_0 < p_1 \le p_2 \le \dots \le p_k \le p_{k+1} = n$ and $b \in
B_{\lambda_1} \otimes \dots \otimes B_{\lambda_n}$ we define
\[
b_{(p_1,\dots,p_k)} = (b_1, b_{\nu_1}, b_2, b_{\nu_2}, \dots,
b_{k+1}, b_{\nu_{k+1}})
\]
where, for $1 \le i \le k+1$, $b_i \in B_\infty$ and $b_{\nu_i}$ is
a highest weight element in $B_{\lambda_{p_{i-1}+1}} \otimes \dots
\otimes B_{\lambda_{p_i}}$ of weight $\nu_i$ as follows.  We define
the $b_i$ and $b_{\nu_i}$ (as well as intermediate elements $a_i$
and $b'_i$) recursively.  First set
\[
a_1 = \hw b,\quad b_1 = \iota^\infty(b).
\]
Now assume that we have defined a highest weight element $a_i$ of
$B_{\lambda_{p_{i-1}+1}} \otimes \dots \otimes B_{\lambda_n}$. If
$i=k+1$, we set $b_{\nu_{k+1}} = a_{k+1}$. Otherwise we have
\[
a_i = b_{\nu_i} \otimes b'_i
\]
for a highest weight element $b_{\nu_i} \in B_{\lambda_{p_{i-1}+1}}
\otimes \dots \otimes B_{\lambda_{p_i}}$ of weight $\nu_i$ and $b'_i
\in B_{\lambda_{p_i + 1}} \otimes \dots \otimes B_{\lambda_n}$ with
$\varepsilon(b'_i) \le \nu_i$.  We then define
\[
a_{i+1} = \hw b'_i,\quad b_{i+1} = \iota^\infty(b'_i).
\]
Thus we have
\[
b = \tilde b_1 (b_{\nu_1} \otimes \tilde b_2 (b_{\nu_2} \otimes
\tilde b_3( \cdots \otimes \tilde b_{k+1} b_{\nu_{k+1}}) \cdots )).
\]

\item For $k \ge 1$ and $0 = p_0 < p_1 \le p_2 \le \dots \le p_k \le p_{k+1} = n$ and $b \in
B_{\lambda_1} \otimes \dots \otimes B_{\lambda_n}$ we define
\[
b^{(p_1,\dots,p_k)} = (b^1, b^{\nu^1}, b^2, b^{\nu^2}, \dots,
b^{k+1}, b^{\nu^{k+1}})
\]
where, for $1 \le i \le k+1$, $b^i$ are the unique elements of
$B_\infty$ and $b^{\nu^i}$ are the unique highest weight elements of
$B_{\lambda_{p_{i-1}+1}} \otimes \dots \otimes B_{\lambda_{p_i}}$ of
weight $\nu^i$ such that
\[
b = \tilde b^1 (b^{\nu^1} \otimes \tilde b^2 b^{\nu^2} \otimes \dots
\otimes \tilde b^k b^{\nu^k}  \otimes \tilde b^{k+1} b^{\nu^{k+1}}).
\]
\end{enumerate}
Note that in the case $k=1$, we have $b^{(p)} = b_{(p)}$ for $1 \le
p \le n$. Also, if $p_i = p_{i+1}$ for some $i$, then we have a
trivial tensor product crystal appearing in the above definitions
and we set $b_{i+1} = b^{i+1} = b_\infty$, $\nu_{i+1} = \nu^{i+1} =
0$ and $b_{\nu_{i+1}} = b^{\nu^{i+1}} = 0$.  In particular, $b^{(n)}
= b_{(n)}$ is always of the form $(b_1, b_{\nu_1}, b_\infty, 0)$
where $b = \tilde b_1 b_{\nu_1}$.
\end{defin}

Whenever we refer to a sequence $(p_1, \dots, p_k)$ as above, we
will always adopt the convention that $p_0=0$ and $p_{k+1}=n$.  If
for some $V \in \mathcal{V}$ we have a flag
\[
0=V^{k+2} \subseteq V^{\nu^{k+1}} \subseteq V^{k+1} \subseteq
V^{\nu^k} \subseteq V^k \subseteq \dots \subseteq V^{\nu^2}
\subseteq V^2 \subseteq V^{\nu^1} \subseteq V^1 = V
\]
of  $I$-graded subspaces and $0 = p_0 < p_1 \le \dots \le p_k \le
p_{k+1} = n$, we say $(x,s,t) \in \mathbf{M}(V,W)$
\emph{$(p_1,\dots,p_k)$-respects} the flag if for all $1 \le i \le
k+1$ we have
\[
x(V^i) \subseteq V^i, \quad x(V^{\nu^i}) \subseteq V^{\nu^i},\quad
s(W^{p_{i-1}+1,p_i}) \subseteq V^{\nu^i}, \quad t(V^i) \subseteq
W^{p_{i-1}+1,n}.
\]
We say the flag is \emph{$(p_1,\dots, p_k)$-respected} by $(x,s,t)$.
In this case, $(x,s,t)$ induces maps
\begin{align*}
x^{V^{\nu^i}/V^{i+1}} &: V^{\nu^i}/V^{i+1} \to
V^{\nu^i}/V^{i+1}, \\
x^{V^i/V^{\nu^i}} &: V^i/V^{\nu^i} \to V^i/V^{\nu^i}
\\
s^{W^{p_{i-1}+1, p_i}, V^{\nu^i}/V^{i+1}} &: W^{p_{i-1}+1, p_i} \to
V^{\nu^i}/V^{i+1}, \text{ and}\\
t^{V^{\nu^i}/V^{i+1},W^{p_{i-1}+1, p_i}} &: V^{\nu^i}/V^{i+1} \to
W^{p_{i-1}+1, p_i}.
\end{align*}

For $b \in B_{\lambda_1} \otimes \dots \otimes B_{\lambda_n}$, let
\[
\tilde T_b = \{(x,s,t) \in \tilde Y_b\ |\ \varepsilon_i (x,s,t) =
\varepsilon_i(Y_b) = \varepsilon_i(b) \ \forall\ i \in I\}.
\]
Then $\tilde T_b$ is a dense subset of $\tilde Y_b$.

For $b_i \in B_\infty$ and highest weight elements $b_{\nu_i} \in
B_{\nu_i}$, $1 \le i \le k+1$, consider the diagram
\begin{multline} \label{eq:Yp-def-diagram}
    \Lambda(\v_1) \times \widetilde{\mathfrak{T}}
    (\v_{\nu_1}; \w^1,\dots,\w^{p_1}) \times \dots \times \Lambda(\v_{k+1}) \times
    \widetilde{\mathfrak{T}}(\v_{\nu_{k+1}}; \w^{p_k+1}, \dots,
    \w^n) \\
    \xleftarrow{\pi_1} \widetilde{\mathfrak{F}}^{(p_1,\dots,p_n)} (\w^1,\dots,\w^n)
    \xrightarrow{\pi_2} \tT,
\end{multline}
where, for $1 \le i \le k+1$, $\v_i$ is the weight of $b_i$, $\nu_i
= \sum_{j=p_{i-1}+1}^{p_i} \lambda_{\w^j} + \v_{\nu_i}$ is the
weight of $b_{\nu_i}$, and
$\widetilde{\mathfrak{F}}^{(p_1,\dots,p_n)} (\w^1,\dots,\w^n)$ is
the variety parameterizing pairs of $(x,s,t) \in \tT$ and flags
\[
0 = V^{k+2} \subseteq V^{\nu_{k+1}} \subseteq V^{k+1} \subseteq
\dots \subseteq V^{\nu_2} \subseteq V^2 \subseteq V^{\nu_1}
\subseteq V^1 = V
\]
$(p_1,\dots,p_k)$-respected by $(x,s,t)$ with dimensions prescribed
by
\[
\dim V^i/V^{\nu_i} = \v_i,\quad \dim V^{\nu_i}/V^{i+1} =
\v_{\nu_i},\quad 1 \le i \le k+1.
\]
The projection $\pi_2$ forgets the flag, while $\pi_1$ is given by
assigning the corresponding induced maps to $(x,s,t)$ and the flag
as above.

\begin{defin} \label{def:Yp}
Let $0 = p_0 < p_1 \le \dots \le p_k \le p_{k+1} = n$.  For $b_i \in
B_\infty$ and highest weight elements $b_{\nu_i} \in B_{\nu_i}$, $1
\le i \le k+1$, let $\pi_1$ and $\pi_2$ be the projections
of~\eqref{eq:Yp-def-diagram}.  Define $\mathcal{Y}^{(p_1, \dots,
p_k)}(b_1, b_{\nu_1}, \dots, b_{k+1}, b_{\nu_{k+1}})$ to be the set
of irreducible components contained in the closure of
\[
    \pi_2 (\pi_1^{-1} (X_{b_1} \times \tilde T_{b_{\nu_1}} \times \dots
    \times X_{b_{k+1}} \times \tilde T_{b_{\nu_{k+1}}})).
\]
\end{defin}
Note that, a priori, $\mathcal{Y}^{(p_1, \dots, p_k)}(b_1,
b_{\nu_1}, \dots, b_{k+1}, b_{\nu_{k+1}})$ may be empty or consist
of several irreducible components.

\begin{lem} \label{lem:Y_b-description-nosplit}
For $b \in B_{\lambda_1} \otimes \dots \otimes B_{\lambda_n}$ with
$b^{(n)} = b_{(n)} = (b_1, b_{\nu_1}, b_\infty, 0)$, the set
$\mathcal{Y}^{(n)}(b_1,b_{\nu_1}, b_\infty, 0)$ consists of the
single irreducible component $Y_b$.
\end{lem}

\begin{proof}
Note that we will always take $V^{\nu_2} = V^2=0$ and so it suffices
to consider the subspace $V^{\nu_1} \subseteq V^1 = V$.  The
condition in Definition~\ref{def:Yp} becomes that $x(V^{\nu_1})
\subseteq V^{\nu_1}$, $s(W) \subseteq V^{\nu^1}$ and
\[
(x^{V^{\nu_1}}, s^{W,V^{\nu_1}}, t^{V^{\nu_1},W}) \in \tilde
T_{b_{\nu_1}},\quad x^{V/V^{\nu_1}} \in X_{b_1}.
\]
Now, for $(x',s',t') \in \tilde T_{b_{\nu_1}}$, we have
$\varepsilon_i(x',s',t')=0$ for all $i \in I$.  Thus for all
$(x,s,t)$ and $V^{\nu_1}$ as above, the smallest $x$-invariant
$I$-graded subspace of $V$ containing $\im s$ is $V^{\nu_1}$ (see
Proposition~\ref{prop:x-nilpotent} and
Lemma~\ref{lem:costable=zero-epsilon}). Therefore, $x$-invariant
subspaces of $V$ containing $\im s$ are in natural one-to-one
correspondence with $x^{V/V^{\nu_1}}$-invariant subspaces of
$V/V^{\nu_1}$.  We now show that
\[
\mathcal{Y}^{(n)}(\kf_{i_l} \cdots \kf_{i_1} b_\infty, b_{\nu_1},
b_\infty, 0) = \{Y_{\kf_{i_l} \cdots \kf_{i_1} b_{\nu_1}}\}
\]
(provided $\kf_{i_l} \cdots \kf_{i_1} b_{\nu_1} \ne 0$) by induction
on $l$.  In the case $l=0$, we take $V^{\nu_1}=V$ and the statement
holds trivially.  Now assume that the result holds for some $l$. For
all $i \in I$, since $\iota^\infty_\lambda$ is $\ke_i$-equivariant,
we have
\[
\varepsilon_i (\kf_{i_l} \cdots \kf_{i_1} Y_{b_{\nu_1}}) =
\varepsilon_i (\kf_{i_l} \cdots \kf_{i_1} b_{\nu_1}) = \varepsilon_i
(\kf_{i_l} \cdots \kf_{i_1} b_\infty) = \varepsilon_i (\kf_{i_l}
\cdots \kf_{i_1} X_{b_\infty}).
\]
Then, upon comparison of the definition of the crystal operators on
$B(\infty)$ and $B(\w^1,\dots,\w^n)$ (see
Sections~\ref{sec:geom-real-Binf} and~\ref{sec:TPQV-crystal}), we
see that
\[
\mathcal{Y}^{(n)}(\kf_{i_{l+1}} \cdots \kf_{i_1} b_\infty,
b_{\nu_1}, b_\infty, 0) = \{Y_{\kf_{i_{l+1}} \cdots \kf_{i_1}
b_{\nu_1}}\}
\]
(provided $\kf_{i_{l+1}} \cdots \kf_{i_1} b_{\nu_1} \ne 0$).
\end{proof}

\begin{prop} \label{prop:Y_b-description}
Let $b \in B_{\lambda_1} \otimes \dots \otimes B_{\lambda_n}$ with
$b^{(p_1, \dots, p_k)} = (b^1, b^{\nu^1}, \dots, b^{k+1},
b^{\nu^{k+1}})$ and $b_{(p_1, \dots, p_k)} = (b_1, b_{\nu_1}, \dots,
b_{k+1}, b_{\nu_{k+1}})$ Then
\[
\mathcal{Y}^{(p_1,\dots,p_k)}(b^1, b^{\nu^1}, \dots, b^{k+1},
b^{\nu^{k+1}}) \text{ and } \mathcal{Y}^{(p_1,\dots,p_k)}(b_1,
b_{\nu_1}, \dots, b_{k+1}, b_{\nu_{k+1}})
\]
each consist of the single irreducible component $Y_b$.
\end{prop}

\begin{proof}
It follows from Definition~\ref{def:bp} that $b_1 = b^1$. We first
prove the case $b^1 = b_\infty$, i.e. $b$ is highest weight.
Consider $\mathcal{Y}^{(p_1,\dots,p_k)}(b^1, b^{\nu^1}, \dots,
b^{k+1}, b^{\nu^{k+1}})$.  Now,
\[
b = b^{\nu^1} \otimes \tilde b^2 b^{\nu^2} \otimes \dots \otimes
\tilde b^k b^{\nu^k} \otimes \tilde b^{k+1} b^{\nu^{k+1}}.
\]
Therefore, by Proposition~\ref{prop:fiber-map-general} and the fact
that the fiber bundle \eqref{eq:p1k-fiber-bundle} induces a crystal
isomorphism, we see that $Y_b$ is the unique irreducible component
such that for all $[x,s,t]$ in a dense subset there is a flag of
$I$-graded spaces
\[
0 = V^{k+2} \subseteq V^{k+1} \subseteq \dots \subseteq V^2
\subseteq V^{\nu^1} = V
\]
with
\[
x(V^i) \subseteq V^i,\quad s(W^{p_{i-1}+1,p_i}) \subseteq V^i,\quad
t(V^i) \subseteq W^{p_{i-1}+1,n},\quad 2 \le i \le k+1,
\]
and
\begin{gather*}
(x^{V^i/V^{i+1}}, s^{W^{p_{i-1}+1,p_i}, V^i/V^{i+1}},
t^{V^i/V^{i+1},W^{p_{i-1}+1,p_i}}) \in \tilde Y_{\tilde b^i
b^{\nu^i}},\quad 2 \le i \le k+1,\\
(x^{V^{\nu^1}/V^2}, s^{W^{1,p_1}, V^{\nu^1}/V^2},
t^{V^{\nu^1}/V^2,W^{1,p_1}}) \in \tilde Y_{b^{\nu^1}}.
\end{gather*}
Then, applying Lemma~\ref{lem:Y_b-description-nosplit} to describe
each $\tilde Y_{\tilde b^i b^{\nu^i}}$ for $2 \le i \le k+1$, we
have that $Y_b$ is the unique irreducible component such that in a
dense subset there exists a flag as above and an $I$-graded subspace
$\bar V^{\nu^i} \subseteq V^i/V^{i+1}$ for $2 \le i \le k+1$ such
that
\[
x^{V^i/V^{i+1}}(\bar V^{\nu^i}) \subseteq \bar V^{\nu^i},\quad
s^{W^{p_{i-1}+1,p_i}, V^i/V^{i+1}}(W^{p_{i-1}+1,p_i}) \subseteq \bar
V^{\nu^i}
\]
and if $V^{\nu^i}$ is the preimage of $\bar V^{\nu^i}$ under the
quotient map $V^i \to V^i/V^{i+1}$ then
\[
(x^{V^{\nu^i}/V^{i+1}}, s^{W^{p_{i-1}+1,p_i}, V^{\nu^i}/V^{i+1}},
t^{V^{\nu^i}/V^{i+1}, W^{p_{i-1}+1,p_i}}) \in \tilde
T_{b^{\nu^i}},\quad x^{V^i/V^{\nu^i}} \in X_{b^i}.
\]
Considering the flag
\[
0 = V^{k+1} \subseteq V^{\nu^k} \subseteq V^k \subseteq \dots
\subseteq V^{\nu^1},
\]
we have the result for $b^1 = b_\infty$.

We now consider $\mathcal{Y}^{(p_1,\dots,p_k)}(b_1, b_{\nu_1},
\dots, b_{k+1}, b_{\nu_{k+1}})$, $b_1=b_\infty$.  We have that
\[
b = b_{\nu_1} \otimes \tilde b_2 ( b_{\nu_2} \otimes \tilde b_3 (
\cdots \otimes \tilde b_{k+1} b_{\nu_{k+1}}) \cdots ).
\]
We prove the result by induction on $k$.  The case $k=0$ is just
Lemma~\ref{lem:Y_b-description-nosplit}.  For $k \ge 1$, by
Proposition~\ref{prop:fiber-map-general} and the fact that the fiber
bundle \eqref{eq:p1k-fiber-bundle} induces a crystal isomorphism, we
see that $Y_b$ is the unique irreducible component such that for all
$[x,s,t]$ in a dense subset there is an $I$-graded subspace $V^2
\subseteq V^{\nu_1}$ such that $x(V^2) \subseteq V^2$,
$s(W^{p_1+1,n}) \subseteq V^2$, $t(V^2) \subseteq W^{p_1+1,n}$, and
\[
(x^{V^{\nu_1}/V^2}, s^{W^{1,p_1}, V^{\nu_1}/V^2}, t^{V^{\nu_1}/V^2,
W^{1,p_1}}) \in \tilde Y_{b_{\nu_1}},\quad (x^{V^2}, s^{W^{p_1+1,n},
V^2}, t^{V^2, W^{p_1+1,n}}) \in \tilde Y_{b'}
\]
where $b' = \tilde b_2 ( b_{\nu_2} \otimes \tilde b_3 ( \cdots
\otimes \tilde b_{k+1} b_{\nu_{k+1}}) \cdots )$.  The result then
follows by the induction hypothesis.

We now prove the general case $b^1=\kf_{i_l} \cdots \kf_{i_1}
b_\infty$, that is
\begin{gather*}
\mathcal{Y}^{(p_1,\dots,p_k)}(\kf_{i_l} \cdots \kf_{i_1} b_\infty,
b^{\nu^1}, \dots, b^{k+1}, b^{\nu^{k+1}}) = \{Y_{\kf_{i_l} \cdots
\kf_{i_1} \hw b}\}, \\
\mathcal{Y}^{(p_1,\dots,p_k)}(\kf_{i_l} \cdots \kf_{i_1} b_\infty,
b_{\nu_1}, \dots, b_{k+1}, b_{\nu_{k+1}}) = \{Y_{\kf_{i_l} \cdots
\kf_{i_1} \hw b}\},
\end{gather*}
for all $l$ (provided $b=\kf_{i_l} \cdots \kf_{i_1} \hw b \ne 0$).
The case $l=0$ is what we have just proved.  The inductive step is
analogous to the one in the proof of
Lemma~\ref{lem:Y_b-description-nosplit} and is therefore omitted.
\end{proof}

We denote the unique element of $\mathcal{Y}^{(p_1,\dots,p_k)}(b_1,
b_{\nu_1}, \dots, b_{k+1}, b_{\nu_{k+1}})$ by
$Y^{(p_1,\dots,p_k)}(b_1, b_{\nu_1}, \dots, b_{k+1},
b_{\nu_{k+1}})$. Note that this is only defined if there exists a $b
\in B_{\lambda_1} \otimes \dots \otimes B_{\lambda_n}$ with
$b^{(p_1, \dots, p_k)}$ or $b_{(p_1, \dots, p_k)}$ equal to $(b_1,
b_{\nu_1}, \dots, b_{k+1}, b_{\nu_{k+1}})$.  From now on, when we
write $Y^{(p_1,\dots,p_k)}(b_1, b_{\nu_1}, \dots , b_{k+1},
b_{\nu_{k+1}})$ we will presuppose the existence of such a $b$.

\begin{cor}
If $b \in B_{\lambda_1} \otimes \dots \otimes B_{\lambda_n}$ with
$b^{(p_1,\dots,p_k)} = (b^1, b^{\nu^1}, \dots, b^{k+1},
b^{\nu^{k+1}})$ and $b_{(p_1,\dots,p_k)} = (b_1, b_{\nu_1}, \dots,
b_{k+1}, b_{\nu_{k+1}})$ then
\[
Y^{(p_1,\dots,p_k)}(b^1, b^{\nu^1}, \dots, b^{k+1}, b^{\nu^{k+1}}) =
Y^{(p_1,\dots,p_k)}(b_1, b_{\nu_1}, \dots, b_{k+1}, b_{\nu_{k+1}}).
\]
\end{cor}
\begin{proof}
This follows immediately from
Proposition~\ref{prop:Y_b-description}.
\end{proof}

We note the difference between the two descriptions of $Y_b$ in
Proposition~\ref{prop:Y_b-description}.  Recall the fact, which
follows easily from the tensor product rule for crystals, that any
highest weight element $b \in B_\lambda \otimes B_\mu$ is of the
form $b_\lambda \otimes b'$ where $b_\lambda$ is the highest weight
element of $B_\lambda$ and $b' \in B_\mu$ with $\varepsilon(b') \le
\lambda$. The first description in
Proposition~\ref{prop:Y_b-description} gives the component $Y_b$ in
terms of the expression of $b$ in the form
\[
b = \tilde b_1 (b_{\nu_1} \otimes \tilde b_2 (b_{\nu_2} \otimes
\tilde b_3 ( \cdots \otimes \tilde b_{k+1} b_{\nu_{k+1}}) \cdots ))
\]
whereas the second describes the same irreducible component in terms
of the expression of $b$ in the form
\[
b = \tilde b^1 ( b^{\nu^1} \otimes \tilde b^2 b^{\nu^2}\otimes \dots
\otimes \tilde b^k b^{\nu^k} \otimes \tilde b^{k+1} b^{\nu^{k+1}}).
\]
These two expressions are obtained from repeatedly applying the
above fact to the different bracketings of the tensor product
\begin{gather*}
(B_{\lambda_1} \otimes \dots \otimes B_{\lambda_{p_1}}) \otimes
((B_{\lambda_{p_1+1}} \otimes \dots \otimes B_{\lambda_{p_2}})
\otimes \dots \otimes (B_{\lambda_{p_k+1}} \otimes \dots \otimes
B_{\lambda_n}) \cdots ),\quad \text{and} \\
( \cdots (B_{\lambda_1} \otimes \dots \otimes B_{\lambda_{p_1}})
\otimes \cdots \otimes (B_{\lambda_{p_{k-1}+1}} \otimes \dots
\otimes B_{\lambda_{p_k}})) \otimes (B_{\lambda_{p_k+1}} \otimes
\dots \otimes B_{\lambda_n})
\end{gather*}
respectively.


\subsection{Action of the commutor on tensor product quiver varieties}

We use the isomorphism $\phi$ of \eqref{eq:phi-isom} to define the
action of the crystal commutor on $B(\w^1, \dots, \w^n)$. In
particular, for $1 \le p \le q < r\le n$ we define
\begin{align*}
\sigma_{p,q,r} &: B(\w^1, \dots, \w^n) \to B(\w^1, \dots, \w^{p-1},
\w^{q+1}, \dots, \w^r, \w^p, \dots, \w^q, \w^{r+1}, \dots, \w^n), \\
\sigma_{p,q,r} &= \phi^{-1} \circ (\id^{\otimes (p-1)} \otimes
\sigma_{B_{\lambda_p} \otimes \dots \otimes B_{\lambda_q},
B_{\lambda_{q+1}} \otimes \dots \otimes B_{\lambda_r}} \otimes
\id^{\otimes (n-r)}) \circ \phi,
\end{align*}
where $\lambda_i = \lambda_{\w^i}$ for $1 \le i \le n$.  When $n=2$,
we write $\sigma : B(\w^1,\w^2) \to B(\w^2, \w^1)$ for
$\sigma_{1,1,2}$.

\begin{prop}
Let $b \in B_{\lambda_{\w^1}} \otimes B_{\lambda_{\w^2}}$ with
$b^{(p)} = b_{(p)} = (b_1, b_{\nu_1}, b_2, b_{\nu_2})$. Then $\sigma
(\mathcal{Y}^{(p)}(b_1,b_{\nu_1},b_2,b_{\nu_2}))$ consists of a
single element and coincides with $\mathcal{Y}^{(p)}(b_1, b_{\nu_2},
b_2^*, b_{\nu_1})$.
\end{prop}

\begin{proof}
We have
\begin{align*}
\sigma ( Y^{(p)}(b_1, b_{\nu_1}, b_2, b_{\nu_2})) &= \phi^{-1}
\sigma_{B_{\lambda_{\w^1}} \otimes B_{\lambda_{\w^2}}} \phi
(Y^{(p)}(b_1, b_{\nu_1}, b_2, b_{\nu_2})) \\
&= \phi^{-1} \sigma_{B_{\lambda_{\w^1}} \otimes B_{\lambda_{\w^2}}}
(\tilde b_1 (b_{\nu_1} \otimes \tilde b_2 b_{\nu_2})) \\
&= \phi^{-1} \tilde b_1 \sigma_{B_{\lambda_{\w^1}} \otimes
B_{\lambda_{\w^2}}} (b_{\nu_1} \otimes \tilde b_2 b_{\nu_2}) \\
&= \tilde b_1 \phi^{-1} (b_{\nu_2} \otimes \tilde b_2^* b_{\nu_1}) \\
&= \phi^{-1} (\tilde b_1 (b_{\nu_2} \otimes \tilde b_2^* b_{\nu_1})) \\
&= Y^{(p)}(b_1, b_{\nu_2}, b_2^*, b_{\nu_1}),
\end{align*}
and the result follows.
\end{proof}

\begin{prop} \label{prop:double-sigma-maps}
Let $b \in B_{\lambda_1} \otimes \dots \otimes B_{\lambda_n}$ with
$b_{(p_1, p_2)} = (b_1, b_{\nu_1}, b_2, b_{\nu_2}, b_3, b_{\nu_3})$
and $b^{(p_1, p_2)} = (b^1, b^{\nu^1}, b^2, b^{\nu^2}, b^3,
b^{\nu^3})$. Then
\begin{gather*}
\sigma_{1,p_1,n} \circ \sigma_{p_1+1,p_2,n} (\mathcal{Y}^{(p_1,
p_2)}(b_1, b_{\nu_1}, b_2, b_{\nu_2}, b_3, b_{\nu_3})) \quad
\text{and} \\
\sigma_{1,p_2,n} \circ \sigma_{1,p_1,p_2}
(\mathcal{Y}^{(p_1, p_2)}(b^1, b^{\nu^1}, b^2, b^{\nu^2}, b^3,
b^{\nu^3}))
\end{gather*}
each consist of a single element and coincide with
\[
\mathcal{Y}^{(p_1, p_2)}(b_1, b_{\nu_3}, b_3^*, b_{\nu_2}, b_2^*,
b_{\nu_1}) \quad \text{and} \quad \mathcal{Y}^{(p_1, p_2)}(b^1,
b^{\nu^3}, (b^3)^*, b^{\nu^2}, (b^2)^*, b^{\nu^1})
\]
respectively.
\end{prop}

\begin{proof}
We have
\begin{align*}
\sigma_{1,p_1,n} &\circ \sigma_{p_1+1,p_2,n} (Y^{(p_1, p_2)}(b_1,
b_{\nu_1}, b_2, b_{\nu_2}, b_3, b_{\nu_3})) \\
&= \phi^{-1} (\sigma_{B_{\lambda_1} \otimes \dots \otimes
B_{\lambda_{p_1}}, B_{\lambda_{p_2+1}} \otimes \dots \otimes
B_{\lambda_n} \otimes B_{\lambda_{p_1+1}} \otimes \dots \otimes
B_{\lambda_{p_2}}}) \\
& \qquad (\id^{\otimes p_1} \otimes \sigma_{B_{\lambda_{p_1+1}}
\otimes \dots \otimes B_{\lambda_{p_2}}, B_{\lambda_{p_2+1}} \otimes
\dots \otimes B_{\lambda_n}}) \phi (Y^{(p_1, p_2)}(b_1, b_{\nu_1},
b_2, b_{\nu_2}, b_3, b_{\nu_3})) \\
&= \phi^{-1} (\sigma_{B_{\lambda_1} \otimes \dots \otimes
B_{\lambda_{p_1}}, B_{\lambda_{p_2+1}} \otimes \dots \otimes
B_{\lambda_n} \otimes B_{\lambda_{p_1}+1} \otimes \dots \otimes
B_{\lambda_{p_2}}}) \\
& \qquad (\id^{\otimes p_1} \otimes \sigma_{B_{\lambda_{p_1}+1}
\otimes \dots \otimes B_{\lambda_{p_2}}, B_{\lambda_{p_2+1}} \otimes
\dots \otimes B_{\lambda_n}}) (\tilde b_1 (b_{\nu^1} \otimes \tilde
b_2 (b_{\nu_2} \otimes \tilde b_3 b_{\nu_3}))) \\
&= \phi^{-1} \sigma_{B_{\lambda_1} \otimes \dots \otimes
B_{\lambda_{p_1}}, B_{\lambda_{p_2+1}} \otimes \dots \otimes
B_{\lambda_n} \otimes B_{\lambda_{p_1}+1} \otimes \dots \otimes
B_{\lambda_{p_2}}} (\tilde b_1 (b_{\nu^1} \otimes \tilde b_2
(b_{\nu_3} \otimes \tilde b_3^* b_{\nu_2}))) \\
&= \phi^{-1} (\tilde b_1 ((b_{\nu_3} \otimes \tilde b_3^*
b_{\nu_2}) \otimes \tilde b_2^* b_{\nu_1})) \\
&= Y^{(p_1, p_2)}(b_1, b_{\nu_3}, b_3^*, b_{\nu_2}, b_2^*,
b_{\nu_1})
\end{align*}
and
\begin{align*}
\sigma_{1,p_2,n} &\circ \sigma_{1,p_1,p_2} (Y^{(p_1, p_2)}(b^1,
b^{\nu^1}, b^2, b^{\nu^2}, b^3, b^{\nu^3})) \\
&= \phi^{-1} (\sigma_{B_{\lambda_{p_1+1}} \otimes \dots \otimes
B_{\lambda_{p_2}} \otimes B_{\lambda_1} \otimes \dots \otimes
B_{\lambda_{p_1}}, B_{\lambda_{p_2+1}} \otimes \dots \otimes
B_{\lambda_n}}) \\
& \qquad (\sigma_{B_{\lambda_1} \otimes \dots \otimes
B_{\lambda_{p_1}}, B_{\lambda_{p_1+1}} \otimes \dots \otimes
B_{\lambda_{p_2}}} \otimes \id^{\otimes (n-p_2)}) \phi (Y^{(p_1,
p_2)}(b^1, b^{\nu^1}, b^2, b^{\nu^2}, b^3, b^{\nu^3})) \\
&= \phi^{-1} (\sigma_{B_{\lambda_{p_1+1}} \otimes \dots \otimes
B_{\lambda_{p_2}} \otimes B_{\lambda_1} \otimes \dots \otimes
B_{\lambda_{p_1}}, B_{\lambda_{p_2+1}} \otimes \dots \otimes
B_{\lambda_n}}) \\
& \qquad (\sigma_{B_{\lambda_1} \otimes \dots \otimes
B_{\lambda_{p_1}}, B_{\lambda_{p_1+1}} \otimes \dots \otimes
B_{\lambda_{p_2}}} \otimes \id^{\otimes (n-p_2)}) (\tilde b^1
(b^{\nu^1} \otimes \tilde b^2 b^{\nu^2} \otimes \tilde b^3 b^{\nu^3})) \\
&= \phi^{-1} (\sigma_{B_{\lambda_1+1} \otimes \dots \otimes
B_{\lambda_{p_2}} \otimes B_{\lambda_1} \otimes \dots \otimes
B_{\lambda_{p_1}}, B_{\lambda_{p_2+1}} \otimes \dots \otimes
B_{\lambda_n}}) (\tilde b^1 ((b^{\nu^2} \otimes \widetilde{(b^2)^*}
b^{\nu^1}) \otimes \tilde b^3 b^{\nu^3}) \\
&= \phi^{-1} (\tilde b^1 (b^{\nu^3} \otimes \widetilde{(b^3)^*}
(b^{\nu^2} \otimes \widetilde{(b^2)^*}
b^{\nu^1}))) \\
&= Y^{(p_1, p_2)}(b^1, b^{\nu^3}, (b^3)^*, b^{\nu^2}, (b^2)^*,
b^{\nu^1}),
\end{align*}
and the result follows.
\end{proof}


\section{The cactus relation}
\label{sec:cactus-relation}

In this section we use the geometric description of the crystal
commutor discussed in Section~\ref{sec:geom-commutor} to show that
the commutor satisfies the cactus relation for arbitrary
simply-laced Kac-Moody algebras.  This result will be extended to
symmetrizable Kac-Moody algebras in
Section~\ref{sec:non-simply-laced}

\subsection{An involution of highest weight irreducible components}

In Section~\ref{sec:geom-kash-involution} we described an involution
on the set of irreducible components of Lusztig quiver varieties
corresponding to Kashiwara's involution. We now discuss a similar
involution on Nakajima quiver varieties. Fix Hermitian forms on $V$
and $W$ such that the form on $W$ is compatible with the
decomposition $W = \bigoplus_{i=1}^n W^i$ (that is, vectors in
different summands are orthogonal).  Consider a point $(x,s,t) \in
\mathbf{M}(V,W)$. Then $(x,s,t)^\dag \stackrel{\text{def}}{=}
(x^\dag, t^\dag, s^\dag) \in \mathbf{M}(V,W)$.  Now
\[
\mu(x,s,t)_i^\dag = \left( \sum_{h \in H,\, \inc(h) = i}
\varepsilon(h) x_h x_{\bar h} + st \right)^\dag = \sum_{h \in H,\,
\inc(h) = i} \varepsilon(h) x_{\bar h}^\dag x_h^\dag + t^\dag s^\dag
= \mu(x^\dag, t^\dag, s^\dag)_i
\]
and so
\[
(x,s,t) \in \mu^{-1}(0) \iff (x,s,t)^\dag \in \mu^{-1}(0).
\]
Recall that $(x,s,t) \in \mu^{-1}(0)$ is a stable point if the only
$I$-graded $x$-invariant subspace of $V$ contained in the kernel of
$t$ is zero.  We say that $(x,s,t) \in \mu^{-1}(0)$ is
\emph{costable} if the only $I$-graded $x$-invariant subspace of $V$
containing the image of $s$ is $V$ itself.  Then it is easy to see
that $(x,s,t)$ is stable if and only if $(x,s,t)^\dag$ is costable.

\begin{lem} \label{lem:costable=zero-epsilon}
For $x$ nilpotent, $(x,s,t)$ is costable if and only if
$\varepsilon_i(x,s,t) = 0$ for all $i \in I$.
\end{lem}
A proof of this lemma appears in \cite[Lemma~2.9.4]{Nak01b}.  We
include a proof for completeness.

\begin{proof}
Suppose that for some $i \in I$, $\varepsilon_i(x,s,t) > 0$.  Then
$\im \tau_i \subsetneq V_i$.  Define $V'_j = V_j$ for $j \ne i$ and
$V'_i = \im \tau_i$.  Then $V'$ is an $x$-invariant proper subspace
of $V$ containing $\im s$.  Therefore $(x,s,t)$ is not costable.

Now suppose that $(x,s,t)$ is not costable.  Then there exists a
proper $I$-graded $x$-invariant subspace $S \subsetneq V$ containing
the image of $s$.  Thus $S^\bot$ is a nonzero $I$-graded
$x^\dag$-invariant subspace of $V$.  Since $x$ (and hence $x^\dag$)
is nilpotent, we can choose a minimal $N$ such that $x_{h_N}^\dag
\cdots x_{h_1}^\dag|_{S^\bot} = 0$ for all directed paths $h_1 h_2
\dots h_N$ in our quiver.  By the minimality of $N$, there exists a
directed path $h_1 h_2 \dots h_{N-1}$ such that $x_{h_{N-1}}^\dag
\cdots x_{h_1}^\dag|_{S^\bot}$ is nonzero. Let $v \in S^{\bot}_i$,
$i = \out(h_{N-1})$, be a nonzero vector in the image of this map.
Now, suppose $h \in H$ with $\inc(h) = i$. By our choice of $N$, $v$
is killed by $x_h^\dag$. Then for all $u \in \out(h)$,
\[
\left< x_h(u), v \right> = \left<u, x_h^\dag (v)\right> = \left<u,
0\right> = 0.
\]
Therefore $v \in (\im x_h)^\bot$ for all $h \in H$ with $\inc(h)=i$.
Furthermore, for all $w \in W$, $\left< s(w), v \right>=0$ since
$s(W) \subseteq S$ and $v \in S^{\bot}$.  Thus $0 \ne v \in (\im
\tau_i)^\bot$. Therefore $\im \tau_i \ne V_i$ and so
$\varepsilon_i(x,s,t) > 0$.
\end{proof}

Let $Y \in B(\w^1, \dots, \w^n)$ be a highest weight element.  Then,
$\varepsilon_i(Y) = 0$ for all $i \in I$.  In other words,
$\varepsilon_i([x,s,t]) = 0$ for all $i \in I$ and $[x,s,t]$ in a
dense subset $U$ of $Y$.  Fix $[x,s,t] \in U$.   Recall that
$(x,s,t) \in \mu^{-1}(0)$ implies $(x,s,t)^\dag \in \mu^{-1}(0)$. By
Proposition~\ref{prop:x-nilpotent} and
Lemma~\ref{lem:costable=zero-epsilon}, $(x,s,t)$ is costable and
thus $(x,s,t)^\dag$ is stable. For $g \in G_V$,
\[
\left( g \cdot (x,s,t) \right)^\dag = (gxg^{-1}, gs, tg^{-1})^\dag =
((g^{-1})^\dag x^\dag g^\dag, (g^{-1})^\dag t^\dag, s^\dag g^\dag) =
(g^{-1})^\dag \cdot (x,s,t)^\dag,
\]
and so
\[
G_V \cdot (x,s,t)^\dag = \left(G_V \cdot (x,s,t)\right)^\dag
\stackrel{\text{def}}{=} \{(x', s', t')^\dag\ |\ (x', s', t') \in
G_V \cdot (x,s,t)\}.
\]
Thus $[x,s,t]^\dag = G_V \cdot (x,s,t)^\dag$ is a well defined point
of $\mathfrak{M}(\w)$.  Let
\[
\lambda'(z) = \id_{W^n} \oplus z \id_{W^{n-1}} \oplus \dots \oplus
z^{n-1} \id_{W^1} = z^{n-1} \lambda(z^{-1}).
\]
Then
\begin{align*}
\lambda'(z) * (x,s,t)^\dag &= \lambda'(z) * (x^\dag, t^\dag, s^\dag) \\
&= (x^\dag, t^\dag \lambda'(z)^{-1}, \lambda'(z) s^\dag) \\
&= (x^\dag, (\lambda'(\bar z)^{-1} t)^\dag, (s \lambda'(\bar
z))^\dag) \\
&= (x, s \lambda'(\bar z), \lambda'(\bar z)^{-1} t)^\dag \\
&= \left( \lambda'(\bar z)^{-1} * (x,s,t) \right)^\dag \\
&= \left( \left(\bar z^{1-n} \lambda(\bar z)\right) * (x,s,t)
\right)^\dag \\
&= \left( \bar z^{n-1} \id_{V} \cdot \lambda(\bar z) * (x,s,t)
\right)^\dag \\
&= z^{1-n} \id_V \cdot (\lambda(\bar z) * (x,s,t))^\dag.
\end{align*}
Therefore
\[
\lambda'(z) * [x,s,t]^\dag = \left( \lambda(\bar z) * [x,s,t]
\right)^\dag
\]
and so
\[
\lim_{z \to 0} \lambda'(z) * \pi([x,s,t]^\dag) = \lim_{z \to 0}
\pi(\lambda'(z) * [x,s,t]^\dag) = \left( \lim_{\bar z \to 0} \pi(
\lambda(\bar z) * [x,s,t]) \right)^\dag = 0.
\]
Thus
\[
U^\dag \stackrel{\text{def}}{=} \{[x,s,t]^\dag\ |\ [x,s,t] \in U\}
\]
is a well defined subset of $\mathfrak{T}(\w^n, \dots, \w^1)$. Since
$\T$ and $\mathfrak{T}(\w^n, \dots, \w^1)$ have the same pure
dimension (they are both lagrangian subvarieties of
$\mathfrak{M}(\w)$ \cite[Prop~3.15]{Nak01}), $U^\dag$ is a dense
subset of some irreducible component of $\mathfrak{T}(\w^n, \dots,
\w^1)$ which we will denote by $Y^\dag$.


\subsection{Proof of the cactus relation}

\begin{prop} \label{prop:double-sigma-equality}
Suppose $Y$ is a highest weight element of the crystal $B(\w^1,
\w^2, \w^3)$. Then
\[
\sigma_{1,1,3} \circ \sigma_{2,2,3} (Y) = \sigma_{1,2,3} \circ
\sigma_{1,1,2} (Y) = Y^\dag.
\]
\end{prop}

\begin{proof}
Choose $b \in B_{\lambda_1} \otimes B_{\lambda_2} \otimes
B_{\lambda_3}$ such that $Y=Y_b$, where $\lambda_i = \lambda_{\w^i}$
for $i = 1,2,3$. Then we have
\begin{gather*}
b_{(1,2)} = (b_\infty, b_{\lambda_1}, b_2, b_{\lambda_2}, b_3,
b_{\lambda_3}),\quad Y=Y^{(1,2)}(b_\infty, b_{\lambda_1}, b_2,
b_{\lambda_2}, b_3, b_{\lambda_3}) \\
\quad b^{(1,2)} = (b_\infty, b_{\lambda_1}, b^2, b_{\lambda_2}, b^3,
b_{\lambda_3}), \quad Y=Y^{(1,2)}(b_\infty, b_{\lambda_1}, b^2,
b_{\lambda_2}, b^3, b_{\lambda_3}),
\end{gather*}
for some $b_2, b_3, b^2, b^3$.  Note that it follows from
Definition~\ref{def:bp} that $\nu^i = \nu_i = \lambda_i$ and
$b^{\lambda_i} = b_{\lambda_i}$ for $i = 1,2,3$.

By Proposition~\ref{prop:double-sigma-maps}, we have
\begin{align*}
\sigma_{1,1,3} \circ \sigma_{2,2,3} (Y) &= Y^{(1,2)}(b_\infty,
b_{\lambda_3}, b_3^*, b_{\lambda_2}, b_2^*, b_{\lambda_1}), \\
\sigma_{1,2,3} \circ \sigma_{1,1,2} (Y) &= Y^{(1,2)}(b_\infty,
b_{\lambda_3}, (b^3)^*, b_{\lambda_2}, (b^2)^*, b_{\lambda_1}).
\end{align*}

Recall that $Y^{(1,2)}(b_\infty, b_{\lambda_1}, b_2, b_{\lambda_2},
b_3, b_{\lambda_3})$ is the unique element of
$\mathcal{Y}^{(1,2)}(b_\infty, b_{\lambda_1}, b_2, b_{\lambda_2},
b_3, b_{\lambda_3})$.  Thus, by Definition~\ref{def:Yp},
$Y^{(1,2)}(b_\infty, b_{\lambda_1}, b_2, b_{\lambda_2}, b_3,
b_{\lambda_3})$ is the unique irreducible component of
$\mathfrak{T}(\w^1, \w^2, \w^3)$ such that for all $[x,s,t]$ in a
dense subset of the component, there exists a flag $0 = V^4
\subseteq V^3 \subseteq V^2 = V^1 = V$ such that
\begin{gather} \label{eq:Vflag-conditions}
x(V^i) \subseteq V^i,\quad s(W^i) \subseteq V^{i+1},\quad t(V^i)
\subseteq W^{i,n},\quad  1 \le i \le 3, \\
\nonumber x^{V/V^3} \in X_{b_2},\quad x^{V^3} \in X_{b_3}.
\end{gather}
Note that $Y_{b_{\lambda_i}} = \{0\}$ and so $V^{\nu^i} = V^{i+1}$
for $i=1,2,3$.  Also $b_1 = b_\infty$ and so $X_{b_1} = \{0\}$.

Similarly, $Y^{(1,2)}(b_\infty, b_{\lambda_3}, b_3^*, b_{\lambda_2},
b_2^*, b_{\lambda_1})$ is the unique irreducible component of
$\mathfrak{T}(\w^3, \w^2, \w^1)$ such that for all $[x',s',t']$ in a
dense subset of the component, there exists a flag $0 = S^4
\subseteq S^3 \subseteq S^2 = S^1 = V$ such that
\begin{gather} \label{eq:Sflag-conditions}
x'(S^i) \subseteq S^i,\quad s'(W^{4-i}) \subseteq S^{i+1},\quad
t'(S^i) \subseteq W^{1,4-i},\quad 1 \le i \le 3,\\
\nonumber (x')^{S/S^3} \in X_{b_3^*},\quad (x')^{S^3} \in X_{b_2^*}.
\end{gather}

For a point $(x,s,t)$ with a flag $0 = V^4 \subseteq V^3 \subseteq
V^2 = V^1 = V$ satisfying \eqref{eq:Vflag-conditions}, set $S^1 =
S^2 = V$, $S^3 = (V^3)^\bot$, and $S^4=0$.  We claim that
$(x',s',t') = (x,s,t)^\dag = (x^\dag, t^\dag, s^\dag)$ satisfies
conditions \eqref{eq:Sflag-conditions}. First we have $x^\dag(S^3) =
x^\dag((V^3)^\bot) \subseteq (V^3)^\bot = S^3$ since $x(V^3)
\subseteq V^3$.  And $x^\dag(S^i) \subseteq S^i$ for $i=1,2,4$
trivially.  The condition $t^\dag(W^3) \subseteq V=S^2$ is also
trivial.  We have $t^\dag(W^2) \subseteq (V^3)^\bot = S^3$ since
$t(V^3) \subseteq W^3$ and $t^\dag(W^1) = 0 = S^4$ since $t(V) =
t(V^2) \subseteq W^{2,3}$.  The condition $s^\dag(S^1) \subseteq W =
W^{1,3}$ holds trivially, $s^\dag(S^2) = s^\dag(V) \subseteq
W^{1,2}$ since $s(W^3) \subseteq V^4=0$, and $s^\dag(S^3) =
s^\dag((V^3)^\bot) \subseteq W^1$ since $s(W^{2,3}) \subseteq V^3$.
The final two conditions then follow from the fact that $x \in X_b$
if and only if $x^\dag \in X_{b^*}$ (see
Section~\ref{sec:geom-kash-involution}). Conversely, $(x,s,t)^\dag$
satisfies \eqref{eq:Sflag-conditions} only if $(x,s,t)$ satisfies
\eqref{eq:Vflag-conditions}. Therefore
\[
\sigma_{1,1,3} \circ \sigma_{2,2,3} (Y) = Y^{(1,2)}(b_\infty,
b_{\lambda_3}, b_3^*, b_{\lambda_2}, b_2^*, b_{\lambda_1}) = Y^\dag.
\]
An analogous argument shows that
\[
\sigma_{1,2,3} \circ \sigma_{1,1,2} (Y) = Y^{(1,2)}(b_\infty,
b_{\lambda_3}, (b^3)^*, b_{\lambda_2}, (b^2)^*, b_{\lambda_1}) =
Y^\dag.
\]
\end{proof}

\begin{cor} \label{cor:sigma-cactus-relation}
We have
\[
\sigma_{1,1,3} \circ \sigma_{2,2,3} = \sigma_{1,2,3} \circ
\sigma_{1,1,2} : B(\w^1, \w^2, \w^3) \to B(\w^3, \w^2, \w^1).
\]
\end{cor}

\begin{proof}
Proposition~\ref{prop:double-sigma-equality} asserts that
$\sigma_{1,1,3} \circ \sigma_{2,2,3} = \sigma_{1,2,3} \circ
\sigma_{1,1,2}$ when restricted to highest weight elements.  The
result then follows from the fact that the maps $\sigma_{p,q,r}$ are
crystal morphisms.
\end{proof}

\begin{theo} \label{thm:cactus-symmetric}
For a Kac-Moody algebra with symmetric Cartan matrix and dominant
integral weights $\lambda_1$, $\lambda_2$, $\lambda_3$,
\[
\sigma_{B_{\lambda_1}, B_{\lambda_3} \otimes B_{\lambda_2}} \circ
\left( \id \otimes \sigma_{B_{\lambda_2}, B_{\lambda_3}} \right) =
\sigma_{B_{\lambda_2} \otimes B_{\lambda_1}, B_{\lambda_3}} \circ
\left( \sigma_{B_{\lambda_1}, B_{\lambda_2}} \otimes \id \right).
\]
That is, the crystal commutor satisfies the cactus relation.
\end{theo}

\begin{proof}
For $i = 1,2,3$, choose $\w^i$ such that $\lambda_i =
\lambda_{\w^i}$. Then we have the crystal isomorphism
\[
\phi : B(\w^1, \w^2, \w^3) \to  B_{\lambda_1} \otimes B_{\lambda_2}
\otimes B_{\lambda_3}
\]
of Proposition~\ref{prop:TPQV-crystal}.  The result then follows
from Corollary~\ref{cor:sigma-cactus-relation} and the definition of
$\sigma_{p,q,r}$.
\end{proof}


\section{Extension to symmetrizable Kac-Moody algebras}
\label{sec:non-simply-laced}

We now extend the results of the previous sections to the more
general setting of symmetrizable Kac-Moody algebras, dropping the
restriction that the Cartan matrix is symmetric.  Our main tool will
be a well-known method for obtaining the Cartan matrices, root
systems, Dynkin diagrams, etc. of non-simply-laced type from the
corresponding objects of simply-laced type via an admissible
automorphism or ``folding'' of a Dynkin diagram.  We refer the
reader to \cite{H04,L93,Sav04b} for details.

\subsection{Admissible automorphisms}
\label{sec:admissible-automs}

Let $(I,E)$ be a graph without loops where $I$ is the set of
vertices and $E$ is the set of edges.  We allow multiple edges
between pairs of vertices. The corresponding symmetric generalized
Cartan matrix is the matrix $A$ indexed by $I$ with entries
\[
a_{ij} = \begin{cases} 2 & i=j \\
-\#\{\text{edges with endpoints $i$ and $j$}\} & i \ne j
\end{cases}.
\]
As usual, let $Q=(I,H)$ be the (double) quiver associated to the
graph.  That is, for each $e \in E$, we have two elements of $H$
arising from the two possible orientations of $e$.  Let $\g(Q)$
denote the symmetric Kac-Moody algebra associated to the above
Cartan matrix, with root system $\Delta(Q)$ (see \cite{K}).

An \emph{admissible automorphism} $\a$ of $Q$ is an automorphism of
the underlying graph such that no edge joins two vertices in the
same $\a$-orbit. Let $\I$ denote the set of vertex $\a$-orbits.
Following \cite{L93} we construct a symmetric matrix $M$ indexed by
$\I$.  The $(\i,\j)$ entry of $M$ is defined to be
\[
m_{\i\j} = \begin{cases} 2\#\{\text{vertices in $\i$th orbit}\} &
\i = \j \\
-\#\{\text{edges joining a vertex in $\i$th orbit and a vertex in
$\j$th orbit}\} & \i \ne \j \end{cases}.
\]
Then let
\[
d_\i = m_{\i \i}/2 = \#\{\text{vertices in $\i$th orbit}\}
\]
and set $D = \diag(d_\i)$.  Then $C = D^{-1}M$ is a symmetrizable
generalized Cartan matrix.  Let $\Gamma$ denote the corresponding
valued graph.  That is, $\Gamma$ has vertex set $\I$ and whenever
$c_{\i \j} \ne 0$, we draw an edge joining $\i$ and $\j$ equipped
with the ordered pair $(|c_{\j \i}|,|c_{\i \j}|)$.  It is known
\cite[Prop 14.1.2]{L93} that any symmetrizable generalized Cartan
matrix (and corresponding valued graph) can be obtained from a pair
$(Q,\a)$ in this way.  The fact that $\a$ is admissible ensures that
$\Gamma$ has no vertex loops.  Let $\g(\Gamma)$ be the Kac-Moody
algebra associated to $C$, with root system $\Delta(\Gamma)$.

Let $(-,-)_Q$ and $(-,-)_\Gamma$ be the symmetric bilinear forms
determined by the matrices $A$ and $M$ respectively.  The
automorphism $\a$ acts naturally on the root lattice $\Z I$ for $Q$,
and $(-,-)_Q$ is $\a$-invariant.  There is a canonical bijection
\[
f: (\Z I)^\a \to \Z \I,\quad f(\beta)_\i = \beta_i \text{ for any }
i \in \i,
\]
from the fixed points in the root lattice for $Q$ to the root
lattice for $\Gamma$.  We will often suppress the bijection $f$ and
consider the root lattice of $\Gamma$ to be the fixed points in the
root lattice for $Q$.  In particular, we have the simple roots for
$\Gamma$ given by
\[
\alpha_\i = \sum_{i \in \i} \alpha_i.
\]
We also define
\[
h_\i = \frac{1}{d_\i} \sum_{i \in \i} h_i.
\]
Then the entries of $C$ are given by $c_{\i \j} = \left< \alpha_\i,
h_\j \right>$.

It was shown in \cite{Kas96} (see also \cite[Lemma~5.1]{Sav04b})
that for vertices $i$ and $j$ in the same orbit $\i$, we have
\[
\ke_i \ke_j = \ke_j \ke_i,\quad \kf_i \kf_j = \kf_j \kf_i
\]
and for any $\g(Q)$-crystal the operators
\[
\ke_{\i} = \prod_{i \in \i} \ke_i,\quad \kf_{\i} = \prod_{i \in \i}
\kf_i
\]
are well defined.  If $B^Q_\infty$ is the $\g(Q)$-crystal
corresponding to the crystal base of $U^-_q(\g(Q))$, then for $b \in
B^Q_\infty$, we also define
\[
\varepsilon_\i(b) = \max \{k \ge 0\ |\ \ke_\i^k b \ne 0\},\quad
\varphi_\i(b) = \varepsilon_\i(b) + \left< h_\i, \wt(b) \right>.
\]
Let $B^\Gamma_\infty$ be the subset of $B^Q_\infty$ generated by the
$\kf_\i$, $\i \in \I$ acting on the highest weight element $b_\infty
\in B^Q_\infty$.  If we restrict the map $\wt : B^Q_\infty \to
P(Q)$, where $P(Q)$ is the weight lattice of $\g(Q)$, to the subset
$B^\Gamma_\infty$, the image lies in the subset of $P(Q)$ invariant
under the natural action of $\a$.  We can therefore view it as a map
$\wt : B^\Gamma_\infty \to P(\Gamma)$ where $P(\Gamma)$ is the
weight lattice of $\g(\Gamma)$.

\begin{prop}
The set $B^\Gamma_\infty$ along with the maps $\ke_\i$, $\kf_\i$,
$\varepsilon_\i$, $\varphi_\i$, $\i \in \I$ and $\wt$ defined above
is a $\g(\Gamma)$-crystal isomorphic to the crystal associated to
the crystal base of $U^-_q(\g(\Gamma))$.
\end{prop}
\begin{proof}
This proposition was proven in \cite{Kas96}.  See also
\cite[Prop~5.2,~Prop~5.5]{Sav04b}.
\end{proof}

Let $\lambda \in P(Q)^+$ be a dominant integral weight of $\g(Q)$
such that $\a(\lambda) = \lambda$.  Thus we can also think of
$\lambda$ as a dominant integral weight of $\g(\Gamma)$.  Let
$B^Q_\lambda$ denote the $\g(Q)$-crystal corresponding to the
irreducible highest weight representation with highest weight
$\lambda$. Let $B^\Gamma_\lambda$ be the subset of $B^Q_\lambda$
generated by the $\kf_\i$, $\i \in \I$, acting on the highest weight
element $b_\lambda$ of $B^Q_\lambda$.  If we restrict the map $\wt :
B^Q_\lambda \to P(Q)$ to the subset $B^\Gamma_\lambda$, then the
image lies in the subset of $P(Q)$ that is invariant under the
action of $\a$.  Thus we can view it as a map $\wt :
B^\Gamma_\lambda \to P(\Gamma)$.

\begin{prop}
The set $B^\Gamma_\lambda$ along with the maps $\ke_\i$, $\kf_\i$,
$\varepsilon_\i$, $\varphi_\i$, $\i \in \I$ and $\wt$ defined above
is a $\g(\Gamma)$-crystal isomorphic to the $\g(\Gamma)$-crystal
corresponding to the irreducible highest weight representation of
$U_q(\g(\Gamma))$ with highest weight $\lambda$.
\end{prop}
\begin{proof}
This proposition was proven in \cite{Kas96}.  See also
\cite[Prop~7.1,~Prop~7.4]{Sav04b}.
\end{proof}


\subsection{The cactus relation for symmetrizable Kac-Moody algebras}

Recall the definition of Kashiwara's involution $*$ in
Section~\ref{sec:kash-involution}.  It is easily seen that when $* :
B^Q_\infty \to B^Q_\infty$ is restricted to $B^\Gamma_\infty
\subseteq B^Q_\infty$, it induces an involution $* : B^\Gamma_\infty
\to B^\Gamma_\infty$ and that this corresponds to Kashiwara's
involution on $B^\Gamma_\infty$, considered as a
$\g(\Gamma)$-crystal.

Let $\lambda, \mu \in P(Q)^+$ be dominant integral weights of
$\g(Q)$ fixed by $\a$.  Thus they can also be viewed as dominant
integral weights of $\g(\Gamma)$.

\begin{lem} \label{lem:tensor-product-a-invariant}
Let $(B_\lambda \otimes B_\mu)^\Gamma$ be the
$\g(\Gamma)$-subcrystal of $B_\lambda \otimes B_\mu$ generated by
the highest weight element $b_\lambda \otimes b_\mu$.  Then
\[
(B_\lambda \otimes B_\mu)^\Gamma = \{b \otimes b'\ |\ b \in
B^\Gamma_\lambda,\ b' \in B^\Gamma_\mu\}.
\]
\end{lem}

\begin{proof}
It suffices to show that for all $\i \in \I$ and $b \otimes b'$ with
$b \in B^\Gamma_\lambda$, $b' \in B^\Gamma_\mu$ and $\kf_\i (b
\otimes b') = b_1 \otimes b'_1$, we have $b_1 \in B^\Gamma_\lambda$
and $b'_1 \in B^\Gamma_\mu$.  Since $b \in B^\Gamma_\lambda$ and $b'
\in B^\Gamma_\mu$, we have
\[
\varepsilon_i(b) = \varepsilon_j(b),\quad \varepsilon_i(b') =
\varepsilon_j(b'),\quad \varphi_i(b) = \varphi_j(b),\quad
\varphi_i(b') = \varphi_j(b'),\quad \text{for } i,j \in \i.
\]
It also follows from the results of
Section~\ref{sec:admissible-automs} that for an element $a$ of
$B^\Gamma_\lambda$ or $B^\Gamma_\mu$ and $i, j \in \i$, $i \ne j$,
\begin{align*}
\varepsilon_i(a) &= \varepsilon_i(\ke_j a) \text{ if } \ke_j a \ne
0,\\
\varepsilon_i(a) &= \varepsilon_i(\kf_j a) \text{ if } \kf_j a \ne
0,\\
\varphi_i(a) &= \varphi_i(\ke_j a) \text{ if } \ke_j a \ne 0,
\text{ and} \\
\varphi_i(a) &= \varphi_i(\kf_j a) \text{ if } \kf_j a \ne 0.
\end{align*}
It follows from the tensor product rule for crystals that
\[
\kf_\i (b \otimes b') = \left(\prod_{i \in \i} \kf_i \right) (b
\otimes b') =
\begin{cases} \left( \left(\prod_{i \in \i} \kf_i \right) b \right)
\otimes b', \text{ or} \\ b \otimes \left(\prod_{i \in \i} \kf_i
\right) b'
\end{cases}
\]
and the result follows.
\end{proof}

It follows from the above results that the crystal commutor
$\sigma_{B^\Gamma_\lambda, B^\Gamma_\mu}$ is obtained by the
restriction of the crystal commutor $\sigma_{B^Q_\lambda, B^Q_\mu}$
to $B^\Gamma_\lambda \otimes B^\Gamma_\mu \subseteq B^Q_\lambda
\otimes B^Q_\mu$.  We thus have the following theorem.

\begin{theo} \label{thm:cactus-symmetrizable}
For a symmetrizable Kac-Moody algebra with dominant integral weights
$\lambda_1$, $\lambda_2$, $\lambda_3$,
\[
\sigma_{B_{\lambda_1}, B_{\lambda_3} \otimes B_{\lambda_2}} \circ
\left( \id \otimes \sigma_{B_{\lambda_2}, B_{\lambda_3}} \right) =
\sigma_{B_{\lambda_2} \otimes B_{\lambda_1}, B_{\lambda_3}} \circ
\left( \sigma_{B_{\lambda_1}, B_{\lambda_2}} \otimes \id \right).
\]
That is, the crystal commutor satisfies the cactus relation.
\end{theo}

\begin{proof}
This follows from Theorem~\ref{thm:cactus-symmetric} and the above
remarks.
\end{proof}


\bibliographystyle{abbrv}
\bibliography{biblist}

\end{document}